\documentclass[11pt,a4paper]{article}
\usepackage{amssymb,latexsym,a4wide,amsmath}
 \DeclareMathOperator{\grad}{grad}

\DeclareMathOperator{\Span}{Span}
\DeclareMathOperator{\Sing}{Sing}
\DeclareMathOperator{\Reg}{Reg}
\DeclareMathOperator{\codim}{codim}
\DeclareMathOperator{\Ker}{Ker}

\DeclareMathOperator{\ord}{ord}
\DeclareMathOperator{\Sym}{Sym}

\def\CC{{\mathbb C}}
\def\RR{{\mathbb R}}
\def\PP{{\mathbb P}}

\def\bull{\vrule height .9ex width .8ex depth -.1ex }

\newtheorem{formula}{}[section]
\newtheorem{proposition}[formula]{Proposition}
\newtheorem{definition}[formula]{Definition}
\newtheorem{corollary}[formula]{Corollary}
\newtheorem{remark}[formula]{Remark}
\newtheorem{lemma}[formula]{Lemma}
\newtheorem{claim}[formula]{Claim}
\newtheorem{theorem}[formula]{Theorem}

\def\thrm{\begin{theorem}}
\def\thrml#1{\begin{theorem}\label{#1}}
\def\ethrm{\end{theorem}}
\def\rmrk{\begin{remark}}
\def\rmrkl#1{\begin{remark}\label{#1}}
\def\ermrk{\end{remark}}
\def\dfntn{\begin{definition}}
\def\dfntnl#1{\begin{definition}\label{#1}}
\def\edfntn{\end{definition}}
\def\nmrt{\begin{enumerate}}
\def\enmrt{\end{enumerate}}

\def\qtn{\begin{equation}}
\def\qtnl#1{\begin{equation}\label{#1}}
\def\eqtn{\end{equation}}
\def\lmm{\begin{lemma}}
\def\lmml#1{\begin{lemma}\label{#1}}
\def\elmm{\end{lemma}}
\def\crllr{\begin{corollary}}
\def\crllrl#1{\begin{corollary}\label{#1}}
\def\ecrllr{\end{corollary}}
\begin{document}
\title{Construction of universal Thom-Whitney-a stratifications, \\
their functoriality and Sard-type Theorem for singular varieties}

\author{
Dima Grigoriev \\[-1pt]
\small IRMAR, Universit\'e de Rennes, \\[-3pt]
\small Beaulieu, 35042, Rennes, France\\[-3pt]
{\tt \small dmitry.grigoryev@univ-rennes1.fr}\\[-3pt]
\small http://perso.univ-rennes1.fr/dmitry.grigoryev
\and
Pierre Milman\\[-1pt]
\small Departement of Mathematics,\\[-3pt]
\small University of Toronto, \\[-3pt]
\small 40 St. George Street, Toronto,\\[-3pt]
\small Ontario M5S 2E4, Canada\\[-3pt]
{\tt \small milman@math.toronto.edu}\\[-3pt]
}
\date{}
\maketitle

\begin{abstract}
{\bf Construction.} For a dominating polynomial (or analytic) 
mapping \mbox{$F: K^n\to K^l$}
with an isolated critical value at $0$ ($K = \RR$ or an algebraically 
closed field of characteristic 
zero) we construct a closed {\it bundle} $G_F \subset T^{*}K^n\ $. 
We restrict $\ G_F \ $ over the critical 
points $\ \Sing(F)\ $ of $\ F\ $ in $\ F^{-1}(0)\ $ and partition $\ \Sing(F)\ $ 
into {\it 'quasistrata'} of points with 
the fibers of $\ G_F\ $ of constant dimension. It turns out that T-W-a 
(Thom and Whitney-a) 
stratifications 'near' $\ F^{-1}(0)\ $ exist iff the fibers of bundle 
$G_F$ are orthogonal to the tangent 
spaces at the smooth points of the quasistrata (e.~g. when $\ l=1\ $). 
Also, the latter are 
the orthogonal complements over an irreducible component $\ S\ $ of 
a quasistratum only if $\ S\ $ is 
{\bf universal} for the class of \mbox{T-W-a} stratifications, meaning 
that for any $\ \{S_j'\}_j\ $ in the class, 
$\ \Sing (F) = \cup _j S'_j \ $, there is a component $\ S'\ $ of an 
$\ S_j'\ $ with $\ S \cap S'$ being open and dense in 
both $\ S\ $ and $\ S'\ $. Construction of {\it Glaeser bundle} $\ G_F\ $ 
involves {\it Glaeser iterations} of 
replacing the fibers of the successive closures by the respective linear 
spans and stabilizes 
after $\ \rho (F) \le 2n\ $ iterations, resulting in $\ \dim (G_F) = n\ $ 
for $\ K \neq \RR \ $.

{\bf Results.} We prove that T-W-a stratifications with only universal 
strata exist iff 
all fibers of $\ G_F\ $ are the orthogonal complements to the respective 
tangent spaces to the 
quasistrata, and then the partition of $\Sing(F)$ by the latter yields 
the coarsest {\it universal 
T-W-a stratification}. (We relax condition of smoothness of strata to 
{\it Gauss regularity}, i.~e. the continuity of their 
Gauss maps.) The proof relies on an extension of a smooth stratum to 
a subvariety with a 
continuous Gauss map and a prescribed tangent bundle over the stratum 
(assuming a version 
of Whitney-a condition). The key ingredient is our version of 
{\bf Sard-type Theorem 
for singular spaces} in which a singular point is considered to be 
noncritical iff nonsingular 
points nearby are 'uniformly noncritical' (e.~g. for a dominating map 
$\ F: X \to Z \ $ 
meaning that the sum of the absolute values of the $l\times l$ minors 
of the Jacobian matrix of 
$\ F \ $, where $\ l = \dim (Z) \ $,
is separated from zero 
by a positive constant). Among examples we include $F:K^5\to K$ that 
does not admit a 
universal T-W-a stratification and a family of 
$\ F_n \colon K^{4n+1} \to K\ $ with $\ \rho (F_n)=n\ $.

{\bf Question.} We wonder whether there can ever be an irreducible 
component of bundle 
$\ G_F\ $ of dimension smaller than $\ n \ $, e.~g. for 
$\ F \colon {\CC}^n\to {\CC} \ $ ?

{\bf Conjecture.} Glaeser bundle $\ G_F\ $ is the intersection of bundles 
associated with
T-W-a stratifications of $\ \Sing(F)\ $ defined over points of $\ \Sing(F)\ $ 
as the orthogonal
complements in the dual of the tangent spaces to the strata passing through 
the respective points.

\end{abstract}

\pagebreak

\tableofcontents

\section*{Introduction}
We consider stratifications of critical points in an isolated critical fiber of 
a dominating polynomial (or analytic) mapping $\ F \colon K^n\to K^l \ $, where 
$K=\RR$ or is an algebraically closed field of characteristic zero, which satisfy 
Thom and Whitney-a conditions. Our main goal is to identify {\bf 'universal strata'}, 
i.~e. such that for every stratification of this type their open and dense subsets 
appear as open dense subsets in appropriate strata of the latter (this gives a hope 
for a solution of the long-standing problem of existence of stratifications with 
double-exponential complexity lower bound). To that end we consider even a larger 
class of Thom-Whitney-a stratifications with the condition of smoothness of strata 
relaxed to {\it Gauss regularity}, i.~e. to a weaker assumption of the existence of 
continuous extensions of their Gauss mappings (sending, by definition, smooth points 
to the tangent spaces at these points) to all points of the strata. Besides being 
{\it Gauss regular} we require strata to be open in their respective closures, 
pairwise disjoint and, of course, to satisfy classical Thom and Whitney-a conditions 
(for the definitions of the latter one may consult for instance \cite{Loo} , 
\cite{Kuo} , \cite{Adam}, \cite{Mac}). {\it Glaeser bundle} $\ G_F \ $ of $\ F \ $ 
is the restriction over the critical points of $\ F \ $ of the subbundle of the 
cotangent bundle which is minimal by inclusion among closed subbundles containing 
the differentials of the component functions of $\ F \ $. 
Construction of $\ G_F \ $ involves {\it Glaeser iterations} of replacing fibers 
of the successive closures by their respective linear spans (see \cite{Glaeser}).

At the first glance it seemed that the Glaeser bundle of the mapping could serve the
purpose of identifying Thom-Whitney-a Gauss regular stratifications
with all strata being universal, namely: by means of partitioning of
the critical locus by dimension of its fibers
(private discussions with A.~Gabrielov, M.~Gromov, M.~Kontsevich, A.~Parusinski and
N.~Vorobjov). But it does not
always work, see example of Subsection~\ref{section7.2}.

Nevertheless, the irreducible subsets (we call them {\it Glaeser components})
over which the fibers of Glaeser bundle are of constant dimension equal
their respective codimension are universal even with respect to the class
of Thom-Whitney-a Gauss regular stratifications, see Corollary~\ref{uni.stratum} .
Thom-Whitney-a stratifications 'near' the critical fiber exist iff the fibers of
Glaeser bundle are orthogonal to the tangent spaces (at the smooth points) of
the {\it quasistrata} of points of constant dimension of fibers of Glaeser
bundle (e.~g. when $\ l=1$, see \cite{Hironaka77}).

Our principal result states that Thom-Whitney-a Gauss regular stratifications
with all strata being universal essentially coincide with the ones derived
from Glaeser bundles by means of the partitioning into the quasistrata described
above. The proof relies on our construction of an extension of a smooth stratum
of a singular locus of a variety to a Gauss regular subvariety with
a prescribed tangent bundle over the stratum under the assumption of
Whitney-a condition on the pair.

To that end our version of a {\bf Sard-type Theorem for singular varieties is
crucial.} Similarly to the classical version its conclusion is that the set of
critical values is 'small', but a singular point is considered to be not critical iff
all 'nearby' nonsingular points are 'uniformly noncritical' (e.~g. for a dominating
map $\ F: X \to Z \ $ meaning that the sum of the absolute values of the $l\times l$
minors of the Jacobian matrix of $\ F \ $, where $\ l = \dim (Z) \ $, not only does not
vanish but, moreover, is separated from zero by a positive constant). Below, following
the setting in which it appears in our paper, we expose a crucial idea of the proof.

Say $X$ is a singular subvariety of an open and dense
${\cal U} \subset \CC^n \ $ and $\ \{ L_j \}_{1 \leq j \leq k}$ is a
collection of functions on $\cal U$ with linearly independent differentials
at each point.
We consider a Sard-type Theorem for the mapping which is the restriction
of the natural projection $X \times \CC^k \rightarrow \CC^k \ $
to $\Lambda_L := \{ (x,c) \in X \times \CC^k : L = 0 \} $, where
$L := \sum _{1\leq j\leq k} c_jL_j \ $. The content of our Sard-type Theorem in
this setting is that for a 'generic' $c \in \CC^k$ not only $d(L_c|_X) \ $, where
$L_c := L|_c \ $, does not vanish at the smooth points of $\ X \ $ in $\{ L_c = 0 \} \ $,
but also that there is a lower estimate (by a positive constant) on the sizes of
$d(L_c|_X)({\alpha}) \ $, where $\ L_c(\alpha) = 0 \ $, points ${\alpha} \in X$
are smooth and are 'nearby' a singular point $ \beta \in X \ $.
We reduce the latter to a problem in a 'nonsingular' setting by means of
an embedded desingularization $\sigma:{\cal N} \rightarrow {\cal U} \ $
of $\ X$ with an additional property that all $L_j \circ \sigma $ become
(locally) monomials and divide each other (for an appropriate ordering).
We apply the standard Sard-type Theorem in this 'nonsingular setting', i.e.
to the restriction of the natural projection $N \times \CC^k \rightarrow \CC^k $
(where smooth $N := \overline{\sigma^{-1}(X \setminus \sigma (Sing(\sigma)))} \subset \cal N$
desingularizes $\ X $) to a smooth hypersurface \mbox{$\Lambda :=
\overline{{\tilde\sigma}^{-1}(\Lambda_L \setminus {\tilde\sigma} (Sing({\tilde\sigma})))} \ $,}
where ${\tilde \sigma} :=$
i\mbox{$\sigma \times id \colon {\cal N} \times \CC^k \rightarrow {\cal U} \times \CC^k \ $.}
Consequently the hypersurfaces
$\Lambda_{(c)} := \Lambda_L \cap (X \times \{ c \})$
are nonsingular off $\sigma (Sing(\sigma))$ for 'generic' $c \in \CC^k \ $.
It remains to carry out the required estimate for ${\alpha} = \sigma(a)$
with noncritical for map $\sigma$ points $\ a \ $ 'nearby' a critical
(also for $\sigma$) point $b \in \{ L_c = 0 \} \ $, where $\beta = \sigma(b) \ $.
(Note that $({\cup}_{1 \leq j \leq k} \Lambda_j) \cup Sing(\sigma) \ $, where
$\Lambda_j := \overline{\sigma^{-1}(\{ L_j = 0 \} \setminus \sigma(Sing(\sigma)))} \ $,
for an appropriate choice of local coordinates is a union of coordinate hyperplanes,
below called 'exceptional'.)

'The crucial idea' can be exposed now as 'an estimate via a logarithmic differentiation':

\noindent we introduce a metric on ${\cal N} \setminus Sing(\sigma)$ 'nearby'
point $\ b \ $ for a choice of local coordinates $x_i \ $, such that the
'exceptional' hyperplanes containing $b$ are $\{ x_i = 0 \}$ for $1 \leq i \leq q$
and one of the remaining $n-q$ coordinates is a local equation of $\Lambda_c \ $,
by 'declaring' collection
$\{ dx_i/x_i \}_{1 \leq i \leq q} \cup \{ dx_j \}_{q+1 \leq j \leq n}$
to be orthonormal; we also introduce a norm on spans of
$\{ dL_j (\sigma (a)) \}_{1 \leq j \leq k}$ by 'declaring' these collections to
be orthonormal. The required estimate follows from the bound (up to a multiplicative
constant) by the size of ${\sigma}^{*}_a ( d(L_c|_X)|_{\sigma(a)})$ on the norms of
the restrictions ${\sigma}^{*}_a :
Span \{ dL_j (\sigma (a)) \}_{1 \leq j \leq k} \rightarrow T^{*}_a ({\cal N})|_{T_a(N)} $
of the pull back by $\sigma$, which at this point is an easy consequence of the
'logarithmic differentiation'.

We provide various examples of mappings that admit universal
Thom-Whitney-a Gauss regular stratifications, but in general the
question of recognition of an individual universal stratum we
address in a forthcoming manuscript: we will show that the universal
strata with respect to Thom-Whitney-a Gauss regular stratifications
are precisely the Glaeser components over which Glaeser bundle is of
dimension $\ n \ $. The latter Glaeser components we refer to as
{\it Lagrangian} since off their singular locus the restriction of
Glaeser bundle over such components is a Lagrangian submanifold of
$T^{*}K^n$ in the natural symplectic structure of the latter.

In abuse of notation we write $\Sing(F)$ for the critical points of $F\ $ in $\ F^{-1}(0)\ $.
We say that  {\it open in its closure} algebraic (or analytic respectively) set $S$
is Gauss regular provided that there is a (unique) continuation
to all of $S$ of the Gauss map from the nonsingular points $\Reg(S)$ of $S\ $, i.~e.
$S\ni x \mapsto T_x(S)\ $, where $T_x(S)$ denotes the tangent space to $S$ at $x\ $.
In abuse of notation we will denote (for a
Gauss regular $S$ and $a \in \Sing(S) := S\setminus \Reg(S)\ $) by $T_a(S)$ the unique limiting
position at $a\ $ of the tangent spaces $T_x(S) \ $ to $S$ for points $x \in \Reg(S)\ $.
We consider Thom-Whitney-a stratifications $\{S_i\}_i$ of the critical points
$\Sing(F)=\cup _i S_i$ with all $S_i$ being Gauss regular (rather than smooth),
open in their respective closures and pairwise disjoint, and such that $\{S_i\}_i$
satisfy Thom and
Whitney-a conditions.
For brevity sake we call them TWG-stratifications and say that $\{S_i\}_i$ is universal
if all irreducible components $S$ of $S_i$ are universal, i.~e. if for any other TWG-stratification
$\{S_j'\}_j$ of $\Sing(F)=\cup _j S_j'$ there exists (a unique) $j$ and an irreducible component $S'$ of
$S_j'$ such that $S\cap S'$ is open and dense in both $S$ and $S'$. Throughout the article by an irreducible
component of a constructible set we mean its intersection with an irreducible component of its closure.

Let quasistrata ${\cal G}_r\subset K^n$ consist of the points of $\Sing (F) \ $ whose
fibers of $G_F$ are vector spaces of dimension $r\ $. Assuming Thom stratification
'near' $F^{-1}(0)$ exists, cf. \cite{Hironaka77} (e.~g. when $\ l=1$), it follows
that $r \ge l$  and that the dimensions of quasistrata ${\cal G}_r$ are less or equal
$n-r$ by virtue of Lemma~\ref{equivalence} below.
Constructed bundle $G_F$ is {\it functorial} with respect to isomorphisms
preserving fibers of $F$ 'near' its critical value $0 \ $ (including with respect to
$C^1$ diffeomorphisms when $K$ is $\CC$ or $\RR$),
see Section~\ref{section2}. Construction of Glaeser bundle $G_F$ involves
iterations (starting with $\{(x,\ \Span\{df_j(x)\}_{1\le j\le l})\}_{x\in K^n} \ $,
where $\Span$ denotes the $K$-linear hull of a family of vectors in $(T_xK^n)^*\ $)
of replacing the fibers of the successive closures by their linear spans and
stabilizes after $\rho (F) \le 2n$ iterations (see \cite{Milman03}), resulting in
$\dim (G_F) = n\ $ for $\ K \neq \RR$ (see Claim~\ref{singular} and Remark~\ref{real}).

The principal purpose of the paper is to provide a constructive criterium of the existence
of a universal TWG-stratification $\{S_i\}_i \ $.
Our main result states that $\Sing(F)$ admits a universal TWG-stratification if and only if
manifolds $\Reg(G_F|_{{\cal G}_r})$ are {\it Lagrangian} in
$K^n\times (K^n)^*$ in the natural symplectic structure of the latter.
Moreover, for universal TWG-stratifications $\{S_i\}_i$ partitions $\{S_{(m)}\}_m$
of $\Sing(F)$ obtained by replacing all $S_i$ of the same dimension $m$ with their union $S_{(m)}$
results in a universal TWG-stratification and coincides with the functorial partition
$\{{\cal G}_r\}_{l \leq r \leq n}$ of $\Sing(F) \ $, which is then the {\it coarsest} among all universal
TWG-stratifications.

A simpler implication that if all $\Reg(G_F|_{{\cal G}_r})$ are Lagrangian then
$\{{\cal G}_k\}_{ l \leq k \leq n}$ is a universal TWG-stratification
we establish in Section~\ref{section3}.
When the latter takes place we would refer to $\{{\cal G}_k\}_{ l \leq k \leq n}$
as a {\it functorial TWG-stratification} (with respect to $F$).

A more difficult converse implication is proved in Sections~\ref{section4}~ and ~\ref{section4.5}.
It relies on Proposition~\ref{continuation} of interest
in its own right. A straightforward generalization of the latter in Theorem~\ref{Continuation} provides
an extension of a (smooth) stratum $\cal G$ of a singular locus of a variety $S$ (algebraic or analytic,
open in its closure and with $\cal G$ being essentially its boundary) to a Gauss regular subvariety
$\cal G^{+}$ of $\overline{S}$ with a prescribed tangent bundle $T_{\cal G}$ over $\cal G \ $ (under
necessary assumptions of our version of Whitney-a condition for the pair of $T_{\cal G}$ over
$\cal G$ and $S$). The key ingredient to both is our version of a Sard-type
Theorem~\ref{sard} for singular varieties. Roughly speaking
Theorem~\ref{sard} asserts that for an irreducible Gauss regular
algebraic (or analytic)  set $S$ its intersection with an appropriate
generic hypersurface (of the same class) is Gauss regular and, more
importantly, the angles between the tangent spaces to $S$ and to the
hypersurface are uniformly separated from 0 on compacts (in a
neighborhood of an open dense subset of any irreducible component of
$\overline{S}\setminus S$).


In Subsection~\ref{section7.index} we construct a family of
$F_n \colon K^{4n+1}\to K$ with the {\it index of stabilization}
$\rho (F_n)\ = n \ $. In Subsection~\ref{section7.2} we prove
that \mbox{$F:=AX^2+2B^2XY+CY^2$}
does not admit a universal TWG-stratification.
Moreover, we show that for an
appropriate variation of the former example an arbitrary hypersurface  appears as
${\cal G}_r$ for some $r$ (see Remark~\ref{ex.sing}). We also consider in
Subsections~\ref{section7.1},~\ref{section7.4} (discriminant-type) examples  for
which $\{{\cal G}_r\}_r$ are functorial TWG-stratifications
(and exhibit these stratifications explicitly).

In abuse of notation in the sequel we identify (occasionaly) the dual $(K^n)^*$ with $K^n$,
the cotangent bundle $T^*(K^n)$  with $K^{2n}$ and also denote $dF(x):=\Span\{\{df_i(x)\}_{1\le i
\le l}\}$. We also denote the variety of zeroes of a polynomial $f$ by $\{f=0\} \ $ and for
the sake of brevity refer to ``Gauss regular'' as ``G-regular''.

\section{Canonical Thom-Whitney-a stratifications}\label{section1}

We recall that in a {\it stratification} $\{S_i\}_i$ of the set
$\Sing(F)=\cup _i S_i$ of critical points of $F$ in $F^{-1}(0)$
(i.~e. the points $x\in F^{-1}(0)$ such that $\dim(dF(x))<l$) each
stratum $S_i$ is assumed to be irreducible (or connected in the
classical euclidean topology for $K = \CC \ $ or $\RR $), open in
its closure and assumed to fulfil the frontier condition: for each
pair $S_i,S_j$ if ${\overline S_i}\cap S_j \neq \emptyset$ then $S_j
\subset {\overline S_i} \ $, as is e.~g. in \cite{Loo}, \cite{Mac}.
Traditionally one assumes each $S_i$ to be smooth.

In the present article for the sake of a concept of universality (and a fortiori functoriality),
i.~e. of a stronger version of canonicity, we relax condition of smoothness and allow $S_i$ to
be \mbox{G-regular.} We consider {\it Gauss regular stratifications}
$Sing(F)=\cup_i S_i$, i.~e. all $S_i$ are G-regular, open
in their respective closures and pairwise disjoint (but neither necessarily irreducible nor
fulfil the frontier condition).
The notions of Thom property with respect to a map $F$ and Whitney-a condition on
stratifications naturally extend to Gauss regular stratifications.

\begin{lemma}\label{criterium}
i) A Thom stratification  exists iff the following condition holds:

(1) any irreducible constructible set $S\subset \Sing(F)$ contains an
open dense subset $S^o\subset \Reg(S)$ such that
if a sequence $\{(x_m \ ,\ dF(x_m))\subset K^{2n}\} _m$ has a limit
$\lim_{m\rightarrow \infty} (x_m \ ,\ dF(x_m))=(x_0 \ ,\ V)$, where $x_0\in S^o$,
$x_m\in K^n \setminus \Sing(F)$ and $V$ is an
\mbox{$l$-dimensional} linear subspace of $(K^n)^*$, then it follows $V \perp T_{x_0}(S^o)$;

$\qquad \quad$ ii) A Thom-Whitney-a stratification  exists iff (1) and  the following
condition hold:

(2) for any smooth irreducible constructible set $M\subset Sing(F)$ and
any irreducible constructible set $S\subset \Sing(F)$ there is an
open dense subset $S^o\subset \Reg(S)$ such that if  a sequence
$\{(x_m \ ,\ V_m)\subset K^n\times (K^n)^*\} _m$
has a limit $\lim _{m\rightarrow \infty} (x_m \ ,\ V_m)=(x_0 \ ,\ V)$, where $x_0\in S^o,\
 x_m\in M$ and subspaces $V_m$ in $(K^n)^*$ are orthogonal to  $T_{x_m}(M)\subset K^n$, then
it follows that subspace $V\subset (K^n)^*$ is orthogonal to  $T_{x_0}(S^o)\subset K^n$.
\end{lemma}

{\bf Proof.} Since the proofs of i) and ii) are similar, we provide only a proof of ii). First
assume that $ \{S_i\}_i$ is a Thom-Whitney-a stratification. Once again the proofs of
properties (1) and (2) are similar and we provide only a proof of (2). Take a unique $S_i$
(respectively, $S_j$) such that $M\cap S_i$ (respectively, $S\cap S_j$) is open and dense
in $M$ (respectively, in $S$). If $S\setminus \overline{S_i}$ is open and dense in $S$ then the choice
of $S^o:=(S_j\cap \Reg(S))\setminus \overline{S_i}$ is as required in (2). On the other hand
the remaining assumptions of (2) can not hold
which makes (2) valid, but vacuous. (Property (1) holds due to the Thom
property of $\{S_i\}_i$.) Otherwise $S\subset \overline{S_i}$ and the choice of $S^o:=S_j\cap
\Reg(S)$ is as required in (1) and in (2) due to the Thom and Whitney-a properties
of  $\{S_i\}_i$ respectively. Indeed, it suffices to replace the sequence of (2) by its
subsequence for which exists $\lim _{m\rightarrow \infty} T_{x_m}(M) =: W \ $, and then
to choose another sequence $\{ x'_m \}_m$ of points in $M\cap S_i$ with the 'distance'
between respective $(x_m \ ,\ T_{x_m}(M)) \ $ and $\ (x'_m \ ,\ T_{x'_m}(M))$ converging to zero.
Then $W = \lim _{m\rightarrow \infty} T_{x'_m}(M)$ and is orthogonal to $V \ $. On the
other hand due to the Whitney-a property of the pair $S_i \ , \ S_j$ it follows that
$W \supset T_{x_0}(S_j) \supset T_{x_0}(S)$ and therefore also $T_{x_0}(S)$ is orthogonal to $V \ $,
as required.

Now we assume that (1) and (2) are valid. We construct  strata
$S_1,S_2,\dots$ by induction on their codimensions, i.~e.
$\codim(S_1)\leq \codim(S_2)\leq \cdots\ $. So assume that
$S_1,\dots,S_k$ are already produced with $\codim(S_k)=r$, set
$\Sing(F)\setminus  (S_1\cup \dots \cup S_k)=:Z$ being of $\codim(Z):=r_1>r$
and that Thom and Whitney-a properties are satisfied for stratification
$\{S_i\}_{1\le i\le k}$ of $\Sing(F)\setminus Z$. Subsequently for every
irreducible component $S$ of $Z$ of $\codim({S})=r_1$ (and by making use
of the noetherian property of the Zariski topology of $S$) we choose a maximal
open subset of $\Reg({{S}})$ which satisfies both property (1) and the property (2)
with respect to the choices of sets $S_i$, for $1\le i\le k$, as the set $M$ of (2).
By additionally choosing each subsequent $S_j$ in $\Sing(F)\setminus (S_1\cup \dots \cup S_{j-1})$
for $k<j\le k_1$ we produce strata $S_{k+1},\dots,S_{k_1}$ of codimensions
$r_1$ with $\codim((\Sing(F)\setminus (S_1\cup \dots \cup S_{k_1}))> r_1$.
Such choice ensures Thom and Whitney-a properties of stratification
$\{S_i\}_{1\le i\le k_1}$ of set $\cup_{1\le i\le k_1} S_i$, as required in the inductive step,
which completes the proof of ii).
\bull

\begin{remark}\label{existence}
It is not true that for $\ l > 1\ $ and $\ 0$ being an isolated critical value of a dominating
polynomial mapping $F: K^n \to K^l$ a stratification that satisfies Thom condition with
respect to $F$ necessarily exists, e.~g~ consider the 'local' blowing up of the origin:
$$F: (z_1,...,z_n) \mapsto (z_1\ ,\ z_1 \cdot z_2\ ,\ \dots\ ,\ z_1 \cdot z_n)\ .$$
\noindent The statement (2) holds, see \cite{Whitney}, \cite{Thom},
\cite{Hironaka73}, \cite{Wall}, \cite{Mac}. For $\ l=1$ statement
(1) holds,
see \cite{Hironaka77}, and for $l > 1\ $ see e.~g.
\cite{Hironaka77}, \cite{Loo}, \cite{Kaloshin} for conditions on $F\
$.
\end{remark}

\begin{remark}\label{similar}
Fix a class of stratifications. A stratification $\{S_i\}_i$ of $\Sing(F)=\cup_i
S_i$ is called {\it canonical} (or {\it minimal}), e.~g. in  \cite{Loo} and  \cite{MR},
if for any other stratification $\{S_i'\}_i\ $ of $\ \Sing(F)=\cup_i S'_i$ in this
class with $\codim(S_1)\leq \codim(S_2)\leq \cdots$ and
$\codim(S'_1)\leq \codim(S'_2)\leq \cdots$ it follows (after possibly
renumbering $\{S'_i\}$) that $S'_1=S_1,\dots,S'_k=S_k$ and
$S'_{k+1}\subsetneq S_{k+1}$. Constructed in the proof of
Lemma~\ref{criterium} Thom and Thom-Whitney-a
stratifications are canonical in the corresponding classes. These respective canonical
stratifications are clearly unique. We extend to Gauss regular stratifications
the concepts and constructions introduced above for stratifications.
\end{remark}

\section{Dual bundles of vector spaces of TWG-stratifications}\label{section2}

In the sequel we will repeatedly apply the following construction.
Let $M \ , \ N$ be constructible sets open in their Zariski closures
(by default we consider Zariski topology, sometimes in the case of
$K$ being $\CC$ or $\RR$ we also use euclidean topology). In the
analytic case we assume alternatively that $M \ , \ N$ are analytic
manifolds. Let $V \ , \ W$ be vector spaces. For a subset ${\cal
T}\subset M\times V$ we denote ${\cal T}^{(0)}={\cal T}$ and by
${\cal T}^{(1)}\subset M\times V$ a bundle of vector spaces whose
fiber ${\cal T}_x^{(1)}$ at a point $x\in M$ is the linear hull of
the fiber $({\overline {\cal T}})_x$ of the closure ${\overline
{\cal T}} \subset M \times V$ \cite{Glaeser}. Defining in a similar
way ${\cal T}^{(p+1)}$ starting with ${\cal T}:={\cal T}^{(p)}\ $,
for $p\geq 0\ ,$ results in an increasing chain of (not necessary
closed) bundles of vector spaces and terminates at ${\cal
T}^{({\rho})}$ such that ${\cal T}^{({\rho})} = {\cal
T}^{({\rho}+1)}$ with ${\rho} \leq 2\dim(V) \ $. We denote $Gl({\cal
T})={\cal T}^{({\rho})}$ and refer to the smallest ${\rho} =
{\rho}({\cal T})$ as the index of stabilization. The so called
'Glaeserization' $Gl({\cal T}) \ $ of $\ {\cal T}$ is the minimal
closed bundle of vector spaces which contains ${\cal T} \ $. We
apply this construction to ${\cal T}=\{(x,dF(x))\}$ where $x$ ranges
over all noncritical points of $F \ $. The result we denote by
$G^{(p)} := G_F^{(p)} := {\cal T}^{(p)}|_{\Sing(F)} \ ,$ for $p\geq
0\ ,$ and $G:= G_F :=Gl({\cal T})|_{\Sing(F)}$ (and still refer to
the smallest ${\rho} = {\rho}(F)$ as the index of stabilization). We
mention that according to \cite{Merle} Thom stratification with
respect to $F$ exists iff $\dim(\overline{G^{(0)}})\le n$, cf.
Remark~\ref{french} and \cite{Hironaka77}. (We do not make use of
the latter criterium in this article.)

Denote $G_x:=\pi ^{-1} (x) \cap G \ $, where
$\pi\colon T^*(K^n)|_{\Sing(F)} \to \Sing(F)$ is the natural projection.
The proofs of the following Proposition and its corollary are straightforward.

\begin{proposition}
Let ${\cal T}_M\subset M\times V,\  {\cal T}_N\subset N\times W $
and $h^{-1} \colon N \to M, H\colon N\times W\rightarrow M\times V$
be  homeomorphisms which commute with the natural projections
$N\times W \to N,\  M\times V\to M\ $. Assume in addition that $H$ is linear on
each fiber of these projections and that $H({\cal T}_N)={\cal T}_M$.
Then $H(Gl({\cal T}_N))=Gl(T_M)$, moreover $H(T_N^{(i)})=T_M^{(i)}$ for every $i$.
\end{proposition}

\begin{corollary}\label{glaeser}
Let $M \ , \ N$ be nonsingular, ${\cal T}_M\subset T^*M,\  {\cal T}_N\subset T^*N$. If
$h:M\rightarrow N$
is an isomorphism such that for the pullback $D^*h$ by $h$ we have $(D^*h)({\cal T}_N)={\cal T}_M$
then $(D^*h)(Gl({\cal T}_N))=Gl({\cal T}_M)$. Moreover,
$(D^*h)({\cal T}_N^{(i)})= {\cal T}_M^{(i)}$ for every $i$.

When $K$ is $\CC$ or $\RR$ it suffices
to assume that $h$ is a $C^1$-diffeomorphism and then constructed bundle $G_F$
and partition $\{{\cal G}_r\}_{l \leq r \leq n}$ of $\Sing(F) \ $ are functorial
with respect to $C^1$ diffeomorphisms preserving fibers of $F$ 'near' its critical value $0 \ $.

\noindent (For an arbitrary $K$ replace ``$C^1$ diffeomorphisms'' above by ``isomorphisms''.)
\end{corollary}

With any Gauss regular stratification ${\cal S}=\{S_i\}_i \ $, where
$Sing(F)=\cup_i S_i \ $, we associate a subbundle $B=B({\cal S})$ of $T^*(K^n)|_{\Sing(F)}$
of vector subspaces of $(K^n)^*$ such that for every $i$ and a smooth point $a\in S_i$
the fiber $B_a:=(T_a(S_i))^{\perp}\subset (K^n)^*$ and for a singular point $a$ of $S_i$
the fiber $B_a$ is defined by continuity, by making use of
 $S_i$ being G-regular. Note that the dimension of fibers $\dim(B_a)=\codim(S_i)$ for $a\in S_i \ $.

\begin{remark}
Note that for any Gauss regular stratification ${\cal S}=\{S_i\}_i$ of $Sing(F)$ bundle
$B({\cal S}) = \cup_i B({\cal S})|_{S_i}$ and for any irreducible component $S$ of an arbitrary
$S_i$ bundle $B({\cal S})|_{S}$ is an irreducible $n$-dimensional
Gauss regular set open in its closure.
\end{remark}

\begin{proposition}\label{almost}
A Gauss regular stratification $\cal S$ satisfies Thom-Whitney-a condition
with respect to $F$ iff $G\subset B$ and $B$ is closed.
\end{proposition}

{\bf Proof.} It follows by a straightforward application of definitions that Thom and
\mbox{Whitney-a} properties for any Gauss regular stratification ${\cal S}=\{S_i\}_i$ of
$Sing(F)$ are equivalent to \mbox{$G^{(1)} \subset B({\cal S})$} and, respectively, that
set $B({\cal S})$ is closed. Due to the definition of bundle $G$ proposition follows. \bull

\begin{corollary}\label{uni.stratum}
It follows due to the preceding Remark and Proposition that all \mbox{$n$-dimensional}
irreducible components of $G$ appear as irreducible components of $B({\cal S})$ for
any \mbox{TWG-stratification} ${\cal S}=\{S_i\}_i$ of $Sing(F)$. Therefore
every irreducible component $\cal G\ $ of $\ {\cal G}_r\ $ with $\ G|_{\cal G}$ being
$n$-dimensional is a universal stratum.
\end{corollary}

Note that $\dim (G) = n$ for $K \neq \RR$ (see Claim~\ref{singular} and Remark~\ref{real}).

\begin{remark}\label{union}
Let $\{S_i\}_i$ be a TWG-stratification of $\Sing(F)$. Then for every $0\leq m\leq n$
the union $\bigcup _{\dim(S_i)=m} S_i$ coincides with
$(\bigcup_{\dim(S_i)\ge m} S_i)\setminus (\bigcup _{\dim(S_i)>m} S_i)$
and therefore is open in its closure. Also due to Proposition~\ref{almost} it is G-regular.
Moreover, if we replace any subfamily of $\{S_i\}_i$ of the same dimension $m$
by its union $S$, we would again obtain a TWG-stratification if
only $S$ is open in its closure.
\end{remark}

\begin{lemma}\label{equivalence}
The following three statements are equivalent:

$\bullet$ a Thom-Whitney-a stratification  exists;

$\bullet$ a TWG-stratification  exists;

$\bullet$ condition (2) of Lemma~\ref{criterium} and the following property hold:

(1') any irreducible constructible set $S\subset \Sing(F)$ contains
an open dense subset $S_0\subset \Reg(S)$ such that for
any $x_0\in S_0$ we have $T_{x_0}(S)\perp G_{x_0}$.
\end{lemma}

{\bf Proof.} For the proof of (1') above note that property (1') with $G_{x_0}$ being replaced
by $G^{(1)}_{x_0}$ is a straightforward consequence of the Thom property of stratification
${\cal S}$ with respect to $F$ and condition (1) of Lemma~\ref{criterium} , which Thom
property implies.
By making use then of condition (2) of Lemma~\ref{criterium} consecutively property (1') with
$G_{x_0}$ being replaced by $G^{(p)}_{x_0} \ $, for $\ p \ge 1 \ $, follows
and implies property (1') as stated, since $G = G^{(p)} \ $ for $\ p = \rho (F) \ $.
Otherwise the proof is similar to that of Lemma~\ref{criterium} with the exception that
we replace $\Reg(S)$ with the maximal (by inclusion) open subset $U$ of $\overline S$
to which by continuity the Gauss map of $S$ uniquely extends from $\Reg(S)$.
\bull
\medskip

Lemma~\ref{criterium} implies (assuming Thom-Whitney-a stratification of $\Sing(F)$ exists) that
$r:=n-\dim(\Sing(F))\ge \min_{a \in \Sing(F)} \{\dim(G_a) \} \ge l$.

\begin{claim}\label{singular}
Assume that Thom stratification of $\Sing(F)$ exists (e.~g. if $l = 1\ $, see
\cite{Hironaka77}), and that $K \neq \RR\ $, then $\Sing(F)=\cup_{j\geq r} {\cal G}_j\ $.
Also, then quasistrata ${\cal G}_j$ are open and dense in irreducible components of
$\ \Sing(F)$ of dimension $n-j$ (if such exist). In particular, appropriate open subsets
of the latter are Lagrangian components of the former with their union being dense
in $\Sing(F) \ $, quasistratum ${\cal G}_r \neq \emptyset$ and $\dim (G_F) = n \ $.
\end{claim}

\begin{remark}\label{real}
In the example of $F: {\RR}^2 \to \RR$ defined by $F := x^3 + x \cdot y^4$ the critical
points $\Sing(F) = \{ 0 \}$, the fiber at $0$ of the Glaeser bundle $G_F$ is spanned by
$dx\ $, i.~e. is \mbox{$1$-dimensional,} and therefore $\dim (G_F) = 1 < 2 =: n \ $.
\end{remark}

{\bf Proof of Claim.} It suffices to verify that a generic point of an irreducible
component of $\ \Sing(F)$ of dimension $n-j$ belongs to ${\cal G}_j \ $, since
the openness is due to
the upper semicontinuity of the function $g\colon x\to \dim(G_x) \ $.

We first reduce to the case of $\ l=1\ $.
Indeed, let $U$ be an open set such that $U \cap \Sing(F)$ is smooth, irreducible
and of dimension $n-j \ $. We may assume w.l.o.g. that $0 \in U \cap \Sing(F)$ and
that for the $1$-st component $f := f_1$ of $F: K^n \to K^l$ the differential $df(0) = 0$
(which anyway holds after a linear coordinate change in the target $K^l$ of map $F$).
By making use of the reduction assumption for $f$ (the case of  $l=1$) it follows that
$(G_f)_a$ are the orthogonal complements of the tangent spaces
$T_a(\Sing(f)) \subset T_a(\Sing(F))\ $ for $\ a$ in an open dense subset $\cal V\ $ of
$\ U \cap \Reg(\Sing(f))\ $. We may also assume by shrinking $U$ and replacing $0\ $, if needed,
that $0 \in \cal V\ $, that $\dim (G_F)_a$ is constant for $a \in U \cap \Sing(F)$ and that
$U \cap \Sing(f) = \cal V$ is smooth, open and dense in an irreducible component of $\Sing(f)$.
Inclusions $\Sing(f) \subset \Sing(F)\ $ and
$\ (G_f)_a \subset (G_F)_a \ $, for $\ a \in \Sing(f)\ $, are straightforward consequences
of the definitions. We continue the proof following

\begin{remark}\label{french}
Note that replacing the assumption of the existence of Thom stratification of
$\Sing (F)$ by the assumption that $\dim (\overline {G^{0}}) \le n$ and following the
proof above would then imply that $(G_F)_a = (G_f)_a\ $, for $\ a \in \cal V\ $, and
moreover that $\dim (U \cap \Sing(F)) = \dim (U \cap \Sing(f))\ $. In particular, it
would follow that $(G_F)_a$ are the orthogonal complements of the tangent spaces
$T_a(\Sing(F)) = T_a(\Sing(f))\ $ for $\ a \in U \cap \Sing(F)\ $, cf. with i) of
Lemma~\ref{criterium} and a criterion $\dim (\overline {G^{0}}) \le n$ for the
existence of Thom stratification of $\Sing (F)$ from \cite{Merle}.
\end{remark}

By making use of the existence of Thom stratification of $\Sing(F)$ and consequently of
(1') of Lemma~\ref{equivalence} applied to $F$ it follows $(G_F)_a$ are orthogonal to
$T_a(\Sing(F))\ $ for $\ a \in U \cap \Sing(f)\ $. Therefore, by making use of the inclusions
above, it follows that $(G_F)_0 = (G_f)_0\ $ and $\ T_0(\Sing(f)) = T_0(\Sing(F))\ $, in
particular implying that $\dim (U \cap \Sing(f)) =$ \mbox{$\dim (U \cap \Sing(F))\ $.} Hence also
$(U \cap \Sing(f)) = (U \cap \Sing(F))\ $, which suffices by making use of the
established above inclusions.

In the case of $l=1$ and by once again making use of (1') of
Lemma~\ref{equivalence} it suffices w.l.o.g. to consider the case of
the restriction  of $\ F$ to a plane of dimension $j$ intersecting
transversally $Z\ $ at $\ a\ $, thus reducing the proof to the case
of $l=1$ and of $a$ being an isolated critical point. In the latter
case it suffices to show that $(G_F)_a = K^n\ $.

If $K$ is algebraically closed our claim follows since for any $c_2\ , \dots \ ,\  c_n \in K$
due to $F_i(a) := {\partial F \over \partial x_i}(a) = 0 \ ,\ 1 \le i \le n\ $, the germ at $a$ of
$\Gamma := \{ F_i - c_i \cdot F_1 = 0\ , 2 \le i \le n \}$ is at least $1$-dimensional, thus producing
$dx_1 + c_2 \cdot dx_2 + \ \cdots \ + c_n \cdot dx_n \ $ in
$\ (\overline{G^{(0)}_F})_a \subset (G_F)_a$ by means of limits of $dF(a)/||dF(a)||$ along
$\Gamma \ $, as required. \bull

\section{Universality and Lagrangian bundles}\label{section3}

Now we introduce a partial order on the class of TWG-stratifications
with respect to $F$ (note that it differs from the order defined in
Ch.1 \cite{Loo}, see Remark~\ref{similar}). For any pair ${\cal
S}=\{S_i\}_i \ , \ {\cal S}'=\{S'_j\}_j \ , \
Sing(F)=\cup_iS_i=\cup_jS_j'$ of TWG-stratifications of $\Sing(F)$
and for every $i$ there exists a unique $j=j(i)$ such that $S_i\cap
S'_j$ is open and dense in $S_i \ $, reciprocately for every $j$
there exists a unique $i=i(j)$ such that  $S_i\cap S'_j$ is open and
dense in $S'_j \ $. We say that $\cal S$ is {\it larger} than ${\cal
S}'$ (i.~e. is 'almost everywhere' finer than $\cal S$) if for every
$i$ equalities $j_0=j(i) \ , \ i=i(j_0)$ hold. Thus universal
TWG-stratification means the largest one.

\begin{proposition}\label{coarser}
For a pair of TWG-stratifications $\cal S$ is larger than  ${\cal S}'$
iff the bundle $B=B({\cal S})\subset B'=B({\cal S}')$.
\end{proposition}

{\bf Proof.} Let $\cal S$ be larger than ${\cal S}'$. For each $i$
we have that $S_i\cap S'_{j_0}$ (where $j_0=j(i)$) is open and dense
in both $S_i,S'_{j_0}$, while $\dim(S_i\cap
S'_{j_0})=\dim(S_i)=\dim(S'_{j_0})$. Therefore, for any point $a\in
S_i\cap S'_{j_0}$ we have $T_a(S_i)=T_a(S'_{j_0}) \ $, i.~e.
$B(S_i)_a=B(S'_{j_0})_a \ $. Hence for any point $b\in S_i$ we obtain
$B_b=B(S_i)_b\subset B'_b$ since the Gauss map of $\overline{S_i}$ is continuous on
$S_i$ and $B'$ is closed due to Proposition~\ref{almost}.

Conversely, let $B\subset B' \ $. For every $S_i$ take $j_0=j(i) \ $, then $S_i\cap S'_{j_0}$
is open and dense in $S_i \ $. It follows that for any point $a\in S_i\cap S'_{j_0}$
inclusion $T_a(S_i)\subset T_a(S'_{j_0})$ holds and therefore $B_a\supset B'_a$ implying that
$B_a=B'_a$ and $\dim(S_i)=\dim(S'_{j_0}) \ $, hence $S_i\cap S'_{j_0}$ is open and dense in
$S'_{j_0} \ $, i.~e. $i(j_0)=i \ $. \bull

Proposition~\ref{coarser} and Remark~\ref{union} imply the following corollary.

\begin{corollary}\label{unique}
i) If for a pair of TWG-stratifications ${\cal S}=\{S_i\}_i$ and
${\cal S}'=\{S'_j\}_j$ (with respect to $F$) equality $B({\cal S})=B({\cal S}')$ holds then the unions
${\cal S}_{(m)}:=\bigcup_{\dim(S_i)=m}S_i=\bigcup_{\dim(S'_j)=m}S'_j$ coincide and are G-regular;

ii) If a universal TWG-stratification ${\cal S}=\{S_i\}_i$
exists then for every $0\leq m\leq n$ the union  ${\cal S}_{(m)}$
is independent of a choice of a universal ~TWG-stratification
and $\{{\cal S}_{(m)}\}_{0\leq m\leq n}$ is
a universal  ~TWG-stratification and is the coarsest universal in the
following sense:
for any universal TWG-stratification ${\cal S}'=\{S'_j\}_j$
and every $0\leq m\leq n$ an equality ${\cal S}_{(m)}={\cal S}'_{(m)}$ holds.
\end{corollary}

For a (constructible) closed subbundle $B\subset T^*(K^n)$ of vector spaces (in the sequel
we shortly call them bundles) we consider its {\it 'quasistrata'},
i.~e.
the constructible sets (open in their respective closures due to the upper-semicontinuity
of the function $\dim_K (B_x)$)
$${\cal B}_k:=\{x\in K^n: \dim_K (B_x)=k\}, \, 0\leq k\leq n.$$
Applying this construction to the bundle $G$ we obtain  quasitrata ${\cal G}_k$.

\begin{remark}
A TWG-stratification exists iff for any point $x \in \Reg({\cal G}_k)$
the fiber $G_x$ is orthogonal to $T_x({\cal G}_k)\ $. Indeed,
the existence of a TWG-stratification implies the desired orthogonality
due to (1') of Lemma 2.7. Conversely, the existence of a TWG-stratification
follows from Lemma 2.7 by making use of the existence of Whitney stratification
\cite{Whitney}, \cite{Kuo}, \cite{Hironaka73}, \cite{Wall}.
\end{remark}

\begin{definition}\label{definition}
We say that irreducible components $\cal B$ of quasistrata ${\cal B}_k \ ,\ 0\leq k\leq n \ , $
are \mbox{\bf Lagrangian} if for points $x\in \Reg({\cal B})$ the tangent spaces $T_x({\cal B})$
are the orthogonal complements of $B_x\ $. We call bundle $B$ Lagrangian if all irreducible
components of ${\cal B}_k \ ,\ 0\leq k\leq n \ , $ are Lagrangian.
\end{definition}

\begin{remark}\label{quasistratum}
For any bundle $B$ Lagrangian components of its quasistrata ${\cal B}_k$ are
G-regular (cf. Remark~\ref{union}) and of dimension $n-k\ $.
\end{remark}

\begin{proposition}\label{up}
If bundle $B$ is Lagrangian then there is a bijective correspondence between the
irreducible components of its quasistrata ${\cal B}_k \ ,\ 0\leq k\leq n \ , $
and the irreducible components of $B$. Also, the irreducible components
$\tilde B$ of $B$ are of dimension $n$ and $\Reg({\tilde B})$ are Lagrangian
submanifolds of $T^*(K^n)$ in the natural symplectic structure of the latter.
\end{proposition}

{\bf Proof}. As a straightforward consequence of Definition~\ref{definition} bundle
$B$ is a union of $n$-dimensional (constructible) sets $B|_{\cal B}$ with
${\cal B}$ being the irreducible components of the quasistrata
${\cal B}_k \ ,\ 0\leq k\leq n \ , $ and $\ \Reg (B|_{\cal B})$ are Lagrangian
submanifolds of $T^*(K^n)\ $. Therefore the closures of $B|_{\cal B}$ are
the irreducible components $\tilde B\ $ of $\ B$ implying the remainder of the
claims of Proposition~\ref{up} as well. \bull

\begin{theorem}\label{lagrangian}
The first two of the following statements are equivalent and imply the third:

(i) bundle $G$ is Lagrangian;

(ii) TWG-stratification of $\Sing(F)$ exists and each irreducible
component of ${\cal G}_k \ $, \mbox{$\ r\leq k\leq n \ $,} is of
dimension $n-k \ $;

(iii) each irreducible component of $G$ is of dimension $n \ $.

\end{theorem}

\begin{remark}
In the example of Remark~\ref{remark.ex} there are only $2$
irreducible components of $G$ and both are of dimension $n=5 \ $,
but $G$ is not Lagrangian.
\end{remark}

{\bf Proof of Theorem~\ref{lagrangian}.} First (i) implies (ii)
since quasistrata $\{{\cal G}_k\}_{r\le k\le n}$ form a
TWG-stratification due to  Proposition~\ref{almost} and
Remark~\ref{quasistratum}. Now assume (ii). Then
(1') of Lemma~\ref{equivalence} implies that for any irreducible
component $\tilde{\cal G}$ of ${\cal G}_k$ there is an open dense
subset ${\tilde{\cal G}}^{(0)}\subset \tilde{\cal G}$ such that
 $T_x({\tilde{\cal G}}) \perp G_x$ holds for any point $x\in {\tilde{\cal G}}^{(0)} \ $.
Since $\dim({\tilde{\cal G}})=n-k$ it follows $G_x$ is the
orthogonal complement to $T_x({\tilde{\cal G}})$ for any point $x\in
{\tilde{\cal G}}^{(0)} \ $, which implies (i). Finally, (i) implies (iii)
is proved in Proposition~\ref{up}. \bull
\smallskip

In the previous section with every TWG-stratification $\cal S$ (with respect
to $F$) we have associated a bundle $B({\cal S})$ such that $B({\cal S})\supset G$ (see
Proposition~\ref{almost}). By construction bundle $B({\cal S})$ is Lagrangian. Conversely, if $B\supset G$
is a Lagrangian bundle then ${\cal S}(B):=\{{\cal B}_k\}_k$ is a
TWG-stratification due to Proposition~\ref{almost} and
Remark~\ref{quasistratum}. We summarize these observations in the following

\begin{theorem}\label{correspondence}
There is a bijective correspondence between TWG-stratifications (with
respect to $F$) and closed Lagrangian subbundles of $T^*(K^n)|_{{\Sing(F)}}$ (which contain $G$).
\end{theorem}

Moreover Propositions~\ref{coarser}, ~\ref{almost},
Theorem~\ref{lagrangian} and Corollary~\ref{unique} imply

\begin{corollary}\label{strongly}
If $G$ is Lagrangian then the corresponding TWG-stratification
$\{{\cal G}_k\}_{r\le k\le n}$ is functorial and is the coarsest
universal.
\end{corollary}

In the next section we establish the converse statement.

\section{A constructive criterium of universality}
\label{section4}

Results of this and of the following section essentially depend on
the validity of the conclusions of Claim~\ref{singular} (which are,
in general, not valid for $K = \RR \ $, cf Remark~\ref{real}). We
therefore additionally assume for the remainder of this article in
the case of $K = \RR$ that bundle $G_F$ is $n$-dimensional over open
dense subsets of every irreducible component of $\Sing (F)\ $. The
latter assumption replaces references below (for $K \neq \RR$) to
Claim~\ref{singular}.

The following Theorem and its Corollary justify the title of the paper.

\begin{theorem}\label{main}
If there exists a universal TWG-stratification
of $\Sing(F)$ then $G$ is Lagrangian.
\end{theorem}

Combining with Corollary~\ref{strongly} it follows

\begin{corollary}
If there exists any universal TWG-stratification of $\Sing(F)$
then $\{{\cal G}_k\}_{r\leq k\leq n}$ is the coarsest universal (and is
functorial).
\end{corollary}

{\bf Proof of Theorem~\ref{main}}. Assume the contrary and let $\cal
G$ be an irreducible component of some ${\cal G}_k \ , \ r\leq k\leq
n$ which is not Lagrangian and with a (lexicographically) maximal
possible pair $(n-k \ ,\ m:=\dim({\cal G})) \ $. We recall (see
Claim~\ref{singular} or in the case $K=\RR$ by an assumption above)
that the minimal $r$ for which ${\cal G}_r\neq \emptyset$ equals
$r=n-\dim(\Sing(F)) \ $. Therefore all irreducible components of
${\cal G}_r$ are Lagrangian since ${\cal G}_r$ is open in $\Sing(F)
\ $, in particular $k>r \ $. We have $m=\dim ({\cal G}) < n-k$ (see
Theorem~\ref{lagrangian}) because condition (1') of
Lemma~\ref{equivalence} implies that $\dim ({\cal G}_t)\leq n-t \ ,
\ r\leq t\leq n \ $. Denote by ${\cal S}=\{S_i\}_i$  a universal
TWG-stratification  of $\Sing(F)=\cup _i S_i$ whose existence is the
assumption  of Theorem~\ref{main} . Below by  an irreducible
component of $\cal S$ we mean an irreducible component of an $S_i \
$.

Let $R\subset \Sing(F) \ $. In the sequel we denote by $G^{\perp}|_R
\subset T(K^n)|_R$ the bundle of vector spaces whose fibers are
the orthogonal complements to the fibers of subbundle $G|_R\subset
T^*(K^n)|_R \ $.

Denote by $W$ the union of all Lagrangian irreducible components of
$\{{\cal G}_t\}_{r\leq t\leq k} \ $. Due to the choice of $\cal G$ we
have $\cup _{r\leq t<k} {\cal G}_t \subset W \ $. On the other hand,
$W$ is the union of all Lagrangian irreducible components of
$\{{\cal G}_t\}_{r\leq t\leq n}$ with dimensions greater or equal to
$n-k \ $. Hence $\dim(\Sing(F)\setminus W)<n-k \ $.

\begin{remark}\label{WandG}
One can produce following the construction in the proof of
Lemma~\ref{criterium} (cf. Remark~\ref{similar}) a TWG-stratification
${\cal S}'=\{S'_j\}_j \ $ of $\ \Sing(F)=\cup_j S_j'$ extending the family of all
irreducible components contained in $W \ $. Then $B(\{S_i\}_i)|_W=G|_W$ due to
Propositions~\ref{almost} and \ref{coarser}. Similarly, $B(\{S_i\}_i)|_{L}=G|_L$
for $L$ being the union (dense in $\Sing(F)$) of all open in $\Sing(F)$
Lagrangian components of appropriate quasistrata ${\cal G}_j \ $ (cf. Claim~\ref{singular}).
\end{remark}

\begin{claim}\label{difficult}
Let $\cal Q$ be an irreducible component of $\cal S \ $. Then either
${\cal Q}\cap W=\emptyset$ or ${\cal Q}$ is an open and dense subset of a Lagrangian
component ${\cal P}\subset W \ $. In particular,
$W$ coincides with the union of an appropriate subfamily of irreducible
components of $\{S_i\}_i \ $.
\end{claim}
{\bf Proof}. Indeed, first consider an
irreducible component $\cal Q$ of $\cal S$ such that ${\cal Q}\cap W$
is dense in $\cal Q$ and denote $t:=n-\dim({\cal Q}) \ $. Since $\cal Q$
is G-regular, $B({\cal S})\supset G$ and $B({\cal S})|_{{\cal Q}\cap W}=G|_{{\cal Q}\cap W}$
it follows that ${\cal Q}\subset \cup_{q\leq t} {\cal G}_q$
and ${\cal Q}\cap W \subset {\cal G}_t$ (in particular $t\le k$).
On the other hand, set ${\cal G}^{(t)}:=\cup_{q\geq t} {\cal G}_q$
is closed (since function $g\colon x\to \dim(G_x)$ is upper semicontinuous)
and therefore ${\cal Q}\subset \overline{{\cal Q}\cap W} \subset {\cal G}^{(t)} \ $.
Hence ${\cal Q}\subset {\cal G}_t \ $.

Consider an irreducible component $\cal P$ of ${\cal G}_t$ such that ${\cal Q}\cap {\cal P}$
is dense in our ${\cal Q} \ $. The latter implies that
$\dim ({\cal P})\geq n-t$ and since ${\cal P}\subset {\cal G}_t$
it follows ( $n-t \geq \dim ({\cal P})$ and therefore) $\dim({\cal P})=n-t \ $.
Thus $\cal P$ is Lagrangian and ${\cal P}\subset W$ (since $t\le k$).
We conclude that
${\cal Q}\subset (\overline{{\cal Q}\cap {\cal P}})\cap {\cal G}_t \subset \overline{\cal P} \cap {\cal G}_t=
{\cal P} \subset W$ and $\dim({\cal Q})=n-t=\dim({\cal P}) \ $, as required.

Now, assume that an irreducible component ${\cal Q}$ of $\cal S$
has a non-empty intersection with a
Lagrangian irreducible component ${\cal P} \subset W$ of ${\cal G}_t$
(and therefore $\dim ({\cal P})=n-t$ for some $t\leq k$). Then,
using $B({\cal S})|_{ {\cal P} \cap  {\cal Q}} = G|_{{\cal P} \cap {\cal Q}}$
and in view of the definition of $B({\cal S}) \ $,
it follows that $\dim({\cal Q})=n-t \ $.
As we have shown above $\dim(\Sing(F)\setminus W)<n-k\le n-t \ $.
Therefore ${\cal Q} \cap W$ is dense in ${\cal Q} \ $.
In the latter case we have already proved that ${\cal Q} \subset W \ $,
which completes the proof of the claim. \bull

\begin{corollary}\label{dense}
Let $\cal Q$ be an irreducible component of $\cal S$ with
$\dim({\cal Q})>\dim({\cal G})$ and $\overline{\cal Q} \supset {\cal G}$
then ${\cal Q} \subset {\cal G}_{n-q} \ $, where $q=\dim({\cal Q})>n-k>\dim({\cal G}) \ $,
and ${\cal Q}\subset W \ $.
\end{corollary}

{\bf Proof.} Due to our assumptions
either ${\cal G}\cap {\cal Q}$
or ${\cal G} \cap (\overline{\cal Q} \setminus {\cal Q})$ is dense
in $\cal G \ $.
If ${\cal Q}\cap W = \emptyset$
then either ${\cal Q}\subset {\cal G}^{(k-1)} \ $ or
$\ {\cal Q}\cap ({\cal G}_k\setminus W)$ is dense in $\cal Q \ $. In the latter
case
$\dim({\cal Q}) \le \dim({\cal G}_k\setminus W)=\dim({\cal G}) \ $, which
is contrary to the choice of $\cal Q \ $. And in the former case
${\cal G}\subset \overline{\cal Q} \subset {\cal G}^{(k-1)}$
contrary to $\cal G$ being an irreducible component of ${\cal G}_k \ $.
Hence ${\cal Q}\cap W \neq \emptyset$ and due to the claim above ${\cal Q}\subset W$. \bull

Consider the union $S^{\cup}$ of all irreducible components $\cal Q \ $
of $\ \cal S$ of the smallest possible dimension $s$
with $\overline{\cal Q}\setminus {\cal Q}$ containing $\cal G \ $.

\begin{remark}\label{Sunion}
Due to the upper semi-continuity of function
$g:x\to \dim(G_x)$ and Claim~\ref{singular}
(or the replacing it assumption when $K = \RR$) the following inclusions hold
${\cal G}\subset \overline{\cup_{r\le t<k} {\cal G}_t} \subset \overline{W}\ $.
Therefore Claim~\ref{difficult}, Corollary~\ref{dense} and Remark~\ref{union}
imply repectively that $S^{\cup}$ is not emply, $S^{\cup}\subset
({\cal G}_{n-s} \cap W) = {\cal G}_{n-s}$ and that $S^{\cup}$ is
G-regular.
\end{remark}

\begin{claim}\label{boundary}
Let $\cal W$ be an irreducible
component of $\overline{S^{\cup}}\setminus S^{\cup}$ such that
$ {\cal W}$ contains $\cal G \ $. Then
${\cal G}$ is dense in $\cal W \ $. (Hence such
$\cal W$ is unique). In particular, $\overline{\cal G}$ is an irreducible component of
$\overline{\overline{S^{\cup}}\setminus S^{\cup}}$ and thus
on an appropriate open neighbourhood $\cal G$ coincides
with its own closure and with $\overline{S^{\cup}}\setminus S^{\cup} \ $.
\end{claim}

{\bf Proof.} Assume the contrary. Then $\dim({\cal W})>\dim({\cal G}) \ $.
Denote by $t_{\cal W}$ the minimal value of $g\colon x\to \dim(G_x) \ $
on $\ \cal W$ (attained on an open dense subset of $\cal W$
in view of the upper semicontinuity of function $g$). Then
$t_{\cal W}\ge t:=n-s=\dim(G_x) \ $ for
$\ x\in S^{\cup}\subset W$ because ${\cal W}\subset (\overline{S^{\cup}}\setminus S^{\cup}) \ $.
Pick an irreducible component ${\cal Q} \ $
of $\ \cal S$ such that ${\cal W}\cap {\cal Q}$ is
dense in $\cal W \ $. Then $\overline{\cal Q} \supset {\cal G}$ and since
$\dim({\cal Q})\geq \dim({\cal W})>\dim({\cal G})$ inclusion ${\cal
Q} \subset W$ holds due to Corollary~\ref{dense}, implying
$( W\cap {\cal G})\supset ({\cal Q}\cap {\cal G}) \ $. Since
${\cal G} \subset ({\cal G}_k \setminus W )$ it follows ${\cal Q}\cap {\cal G}$ is empty, i.~e.
${\cal G}\subset (\overline{\cal Q}\setminus {\cal Q}) \ $. Since also
${\cal Q} \subset W$ and due to the choice of $s$ we conclude that $\dim({\cal Q})\ge s \ $.
On the other hand $n-\dim({\cal Q})=\dim(G_x)=t_{\cal W} \ $ for $\ x\in ({\cal W}\cap {\cal Q})$
by making use of Remark~\ref{WandG} and Claim~\ref{difficult}, which
implies $s=n-t\ge n-t_{\cal W}=\dim({\cal Q}) \ $. Therefore $s=\dim({\cal Q})$ and
both ${\cal Q}\subset S^{\cup}$ and, due to
$\overline{\cal Q}\cap \overline{\cal W}\neq \emptyset \ $,
inequality ${\cal Q} \cap (\overline{S^{\cup}}\setminus S^{\cup}) \neq \emptyset$ holds,
leading to a contradiction. \bull

\begin{corollary}\label{germ}
Let $\cal Q$ be an irreducible component of $\cal S$ of $\dim({\cal Q})=s$
with $\overline{\cal Q} \setminus {\cal Q}\supset \cal G \ $. Let
$S_*:=\overline{\cal Q}\cap S^{\cup}\supset \cal Q \ $.
Then $S_*$ is an irreducible subset of  $W\cap {\cal G}_{n-s}={\cal G}_{n-s} \ $ and
$\ \overline{S_*}\setminus S_*= \cal G=\overline{\cal G}$ in an open neighbourghood
$U_{\cal G} \ $.
\end{corollary}

{\bf Proof.} Inclusion $S_*\subset S^{\cup} \subset W\cap {\cal G}_{n-s}={\cal G}_{n-s}$ is the main
content of Corollary~\ref{dense}.
Note that $S_*$ is irreducible since $\overline{S_*}=\overline{\cal Q} \supset \cal G$ and
that sets ${\cal G}\cap S^{\cup}$ and
$(\overline{S_*}\setminus S_*)\cap S^{\cup}$ are both empty.
Therefore $S_*\cap {\cal G}=\emptyset$ and $(\overline{S^{\cup}}\setminus S^{\cup}) \supset
(\overline{S_*}\setminus S_*) \supset \cal G \ $. Hence due to Claim~\ref{boundary} also
$\overline{S_*}\setminus S_*$ coincides with $\cal G$ on an open neighbourhood of an open
dense subset of $\cal G \ $. \bull

\begin{remark}\label{open}
We may choose an open neighbourhood $U_{\cal G}$ of $\cal G$ so that
${\cal G} \cap U_{\cal G}=\overline{\cal G}\cap U_{\cal G} \ $.
Since
$\overline{\cal Q}\cap U_{\cal G}\supset
{\cal G} \cap U_{\cal G}\neq \emptyset$ it follows that
${\cal Q}\cap U_{\cal G} \neq \emptyset \ $. Consider
$S:=S_*\cap U_{\cal G}\supset {\cal Q}\cap U_{\cal G}$ (as in Corollary~\ref{germ}).
Then $\overline{\cal Q}\supset \overline{S}\supset
 \overline{{\cal Q}\cap U_{\cal G}} =\overline{\cal Q}=\overline{S_*}\ $ (due to
$\cal Q$ being irreducible) and therefore  $\overline{S}=\overline{S_*} \ $ and $\ S$
is irreducible. Hence
${\cal G} \cap U_{\cal G}=
(\overline{S_*}\setminus S_*)\cap U_{\cal G} \supset  (\overline{S}\setminus S)\cap U_{\cal G}
\supset {\cal G} \cap U_{\cal G}\ $,
which implies
\begin{eqnarray}\label{1.2}
(\overline{S}\setminus S)\cap U_{\cal G}= {\cal G} \cap U_{\cal G} =\overline{\cal G}\cap U_{\cal G}
\end{eqnarray}
\noindent and that $S$ is open in its closure. Finally, $S$ is
G-regular (and is a dense subset of a Lagrangian component of ${\cal
G}_{n-s}$) since $S\subset W\cap {\cal G}_{n-s}={\cal G}_{n-s} \ $.
\end{remark}

In the remainder of this and in the following Section we use notation
$\cal G \ $ for $\ {\cal G}\cap U_{\cal G} \ $ and $\ S \ $
for $\ S\cap U_{\cal G}$ from Remark~\ref{open}.

\begin{proposition}\label{continuation}
There is an irreducible G-regular constructible set ${\cal G}^+$
open in its closure such that ${\cal G}^+ \subset \overline S \ $, $\ \dim ({\cal G}^+) =n-k \ $
and $\ {\cal G}^+$ contains an open dense subset of ${\cal G} \ $. Finally
$$G^{\perp}|_{{\cal G}^+ \cap {\cal G}}=T({\cal G}^+)|_{{\cal G}^+ \cap {\cal G}} \ .$$
\end{proposition}

{\bf Deduction of Theorem~\ref{main} from
Proposition~\ref{continuation}.} The bundle of vector spaces associated
(as in  Section~\ref{section2}) with a family
$$W_1=\bigcup_{{\cal Q}\subset W} ({\cal Q}\setminus \overline{{\cal G}^+}) \cup \{{\cal G}^+\}$$
\noindent (where the union ranges, as above, over all irreducible
components $\cal Q \ $ of $\ \cal S$ such that ${\cal Q}\subset W$)
coincides over $W_1 \setminus {\cal G}^+$ with G, is Lagrangian and
is closed due to the latter and Proposition~\ref{continuation}.
Since $W\setminus W_1 \subset \overline{{\cal G}^+}\setminus {\cal
G}^+$ and dimensions of $(\overline{{\cal G}^+} \setminus {\cal
G}^+) \ $ and $\ (\Sing(F)\setminus W)$ are less than $n-k$ it
follows that $\dim (\Sing(F)\setminus W_1)< n-k \ $. Therefore, as
in the Remark~\ref{WandG}, the latter family extends to a
TWG-stratification $\{{\tilde S}_j\}_j \ $ of $\ \Sing(F)=\cup _j {\tilde S}_j \ $.

As we have established above in Claim~\ref{difficult} set $W$
and therefore $Sing(F)\setminus W$ are the unions of several
irreducible components of $\cal S \ $. Hence there exists an
irreducible component $\cal P \ $ of $\ \cal S$ such that
$(\Sing(F)\setminus W)\supset {\cal P} \ $ and
$\ {\cal G} \cap \cal P$ is open and dense in $\cal G \ $. Since
being universal TWG-stratification $\{S_i\}_i$ is larger
than $\{{\tilde S_j}\}_j$ it follows by Proposition~\ref{coarser}  that
for any point $x\in  {\cal G} \cap {\cal G}^+ \cap \cal P$ there is
an inclusion  $B({\cal P})_x\subset B({\cal G}^+)_x=G_x$ for the fibers of $G \ $;
hence $\dim(B({\cal P})_x) \leq \dim (G_x)=k \ $ and $\ \dim({\cal P})\geq n-k \ $. But
on the other hand $\dim({\cal P})\leq \dim((\Sing(F)\setminus W)<n-k \ $.
Thus the assumption (on the first lines of the proof of Theorem~\ref{main})
of the existence of a non Lagrangian component $\cal G \ $ in
$\ \{{\cal G}_j\}_j$ leads to a contradiction, i.~e. $G$ is Lagrangian.
\bull

\section{Sard-type Theorem for singular varieties}\label{section4.5}

Proof of the more difficult implication of our main result Theorem~\ref{main}
we complete in this section. To that end
we prove here Proposition~\ref{continuation}, which essentially
provides an extension of a (smooth) singular locus of an algebraic variety
to a Gauss regular subvariety with a prescribed tangent bundle over singularities.
The main ingredient is our Sard-type Theorem for singular varieties.
\smallskip

To begin with we introduce a generalization of Whitney-a property
for a pair ${\cal G} \ , \ S$ of smooth irreducible algebraic (or analytic
respectively) sets closed in a nonsingular ambient variety $U_{\cal G} \ $, and in
$U_{\cal G} \setminus {\cal G}$ respectively, with ${\cal G}$ being the boundary
of $S$ in $U_{\cal G} \ $. Our generalization requires additional data of a
subbundle $T_{\cal G} \ $ over $\ {\cal G}$ of the tangent bundle $T(U_{\cal G})|_{\cal G} \ $
of $\ U_{\cal G}$ (restricted over ${\cal G}$)
that contains the tangent bundle of ${\cal G} \ $. (To apply the notion in the
setting of Proposition~\ref{continuation} we allow $S$ to be Gauss regular.)
Then our generalized Whitney-a condition is as follows:

W-a) if  a sequence $\{ (x_i, T_{x_i}(S)) \subset S\times T(U_{\cal G})|_S \} _i$
has a limit $\lim _{i\rightarrow \infty} (x_i, T_{x_i}(S)) = (x_0,V) \ $, where
$x_0\in {\cal G}$ and subspace $V\subset T_{x_0}(U_{\cal G}) \ $, then it
follows that subspace $V \supset (T_{\cal G})_{x_0} \ $.
\smallskip

\begin{theorem}\label{Continuation}
Assume ${\cal G} \ , \ S \ , \ U_{\cal G}$ and $T_{\cal G} \subset T(U_{\cal G})|_{\cal G}$
are as in the preceding paragraph and satisfy generalized Whitney-a condition W-a).
Then there is an irreducible Gauss regular closed subvariety ${\cal G}^+$ of
$\overline {S}$ in an open subset $U'_{\cal G}$ of $U_{\cal G}$ that contains
an open dense subset ${\cal G} \cap U'_{\cal G} \ $ of $\ {\cal G}$ and such that
$$T_{\cal G}|_{{\cal G}^+ \cap {\cal G}} = T({\cal G}^+)|_{{\cal G}^+ \cap {\cal G}} \ .$$
\end{theorem}
\smallskip

\begin{remark}\label{specialcase}
Theorem~\ref{Continuation} is a straightforward generalization of Proposition~\ref{continuation}
and a straightforward extension of the proof of the latter below applies to the former.
\end{remark}
\smallskip

{\bf Proof of Proposition~\ref{continuation}.} Throughout the proof
of the Proposition we assume that the field $K=\CC \ $ (or $\ \RR \
$), and afterwards extend the proposition to an arbitrary
algebraically closed field employing the Tarski-Lefschetz principle.

First we construct a $(k+m)\times n \ $ matrix $\ M=(M_{j,i})_{1\leq
j\leq k+m,1\leq i\leq n} \ $ with the entries being polynomials over
$K = \CC \ $ (or $\ \RR$) in $n$ variables such that for a suitable open subset
$V\subset \cal G$ we have
\begin{eqnarray}\label{15}
G^{\perp}|_V=T({\cal G})|_V \oplus \Ker(M)|_V \ .
\end{eqnarray}
\noindent
In particular, the rank of $M$ equals $k+m$ at all  points of $V \ $.

Consider a Noether normalisation $\pi:{\cal G} \rightarrow {K}^m$
being a restriction of a linear projection $\pi:{K}^n\rightarrow
{K}^m \ $. Assuming that ${K}^m\subset {K}^n$, one can represent
${K}^n={K}^m\oplus {K}^{n-m}$ with ${K}^{n-m}=\Ker(\pi)$ and
${K}^m=\pi({K}^n) \ $. We may  assume w.l.o.g. that the first $m$
coordinates are the coordinates of the first summand and the last
$n-m$ coordinates are the coordinates of the second summand.
We choose in the tangent space to ${K}^n$ the respective to these
\mbox{$X$-coordinates} a basis of $\partial \over {\partial X_i} \ $. In abuse of
notation we denote ${K}^{n-m} = T_x({K}^{n-m}) \subset T_x({K}^n)$
for \mbox{points $x \in {K}^{n-m} \ $.}

Take an open subset ${\cal U}\subset {K}^m$ such that (\ref{15}) holds for
$V:=\pi^{-1}({\cal U}) \cap \cal G$ and the dimension of any fiber of the
bundle $$G^{\perp}|_V \cap (V\times {K}^{n-m})$$
\noindent equals $n-k-m \ $, e.~g. any open $\cal U$ such that over $V$ the tangent
spaces to $\cal G$ are mapped onto ${K}^m$ isomorphically would do. Note that since
${\cal G}={\cal G}\cap U_{\cal G}$ it follows that $\pi(U_{\cal G}\cap V)={\cal U} \ $.
Then there is a $(k+m)\times n$ matrix $M$ such that
$$\Ker(M)|_V=G^{\perp}|_V \cap (V\times {K}^{n-m})\ .$$
\noindent Of course we may assume w.l.o.g. that
$M_{j,i}=\delta_{j,i} \ $ for $\ 1\leq j\leq m\ , \ 1\leq i\leq n$ (where $\delta$ denotes
the Kronecker's symbol).
This provides a required matrix $M$, a set $V$ and (\ref{15}).

One can construct (by means of an
interpolation in ${K}^{n-m}$ parametrized by points in ${\cal U}' \ $,
see Appendix) rational in the first $m$ (and polynomial in the last
$n-m$) coordinates functions $L_j(X)\ , \ 1\leq j\leq k\ $,
and an open subset ${\cal U}'\subset \cal U$ such that
 all $L_j\ , \ 1\leq j\leq k \ $, vanish on $V':=\pi^{-1}({\cal U}') \cap \cal G$
(while their denominators do not)
and for every point $x\in V'$ equalities
$${\partial L_j \over \partial X_i}(x)=M_{j+m,i}(x)\ ,\ 1\leq j\leq k\ ,\ m+1\leq i\leq n\ ,$$
\noindent hold.
Multiplying by the common denominator and keeping the same notation for  polynomials
$L_j \ , \ 1\leq j\leq k \ $ we conclude that all $L_j$ vanish on $\cal G \ $,
their differentials $dL_j(x) \ , \ 1\leq j\leq k$ are linearly independent
for any $x\in V'$ and due to (\ref{15})
\begin{eqnarray}\label{0}
\bigcap_{1\leq j\leq k} \Ker(dL_j)|_{V'}\ =\ G^{\perp}|_{V'}\ .
\end{eqnarray}
\noindent Therefore by shrinking neighbourhood $U_{\cal G}$ if necessary
we may assume w.l.o.g.
that $U_{\cal G}\subset \pi^{-1}({\cal U}')$ and
that differentials $dL_1,\dots,dL_k$ are linearly independent at every point in
$U_{\cal G} \ $.

A collection of varieties forms a  {\it normal crossings} at a point
$a$ provided that in appropriate analytic local coordinates centered
at this point every variety from this collection and passing through
$a$  is a coordinate subspace. Of course this property is open with
respect to the choice of points $a$. Due to our choice above, collection of
hypersurfaces $H_j:= \{L_j=0\}\cap U_{\cal G}\ ,\ 1\le j\le k\ ,$
forms normal crossings in $U_{\cal G}\ $, i.~e. at every point of
$U_{\cal G}\ $. Moreover, since $S$ is irreducible (see
Remark~\ref{open}) it follows that the set $\Reg_*(\overline{S})$ of
points of $\overline{S}\cap U_{\cal G}$ at which collection of
$\{H_j\}_{1\le j\le k}$ with $\overline{S}$ forms normal crossings is an open
and dense subset of $\Reg(\overline{S}\cap U_{\cal G})\ $ (since
$\Reg_*(\overline{S})\supset \Reg(S)\setminus \bigcup
_{H_j\not\supset S} H_j \neq \emptyset$). In the sequel we denote
$\Sing_*(\overline{S}):= \overline{S}\cap  U_{\cal G}\setminus
\Reg_*(\overline{S})$.

To complete the proof of our Proposition we will need a Sard-type
Theorem  for {\it  singular varieties}. We observe that due to
Proposition~\ref{almost} and $S$ being dense in a Lagrangian component
of ${\cal G}_{n-s}$ (see Remark~\ref{open}) inclusion
$$\overline{T(S)^{\perp}}|_{\cal G}\subset G|_{\cal G}$$
\noindent holds.
In a version of Sard-type Theorem
below assuming the latter inclusion and (\ref{1.2}) we construct in $S$
a codimension one G-regular subvariety
${\hat S}_{-1}:={\hat S}_{-1}(S)\subset S$ with
$(\overline{{\hat S}_{-1}} \cap S) = {\hat S}_{-1} \ $,
such that $\overline{{\hat S}_{-1}} \supset \cal G$ and
inclusion
$$\overline{T({\hat S}_{-1})^{\perp}}|_{\cal G}\subset {G}|_{\cal G}$$
\noindent holds (thus, the pair ${\hat S}_{-1} \ , \ {\cal G}$ behaves
similarly to the pair
${\hat S}_{-0}:=S \ , \cal G\ $, cf.  items iii)-vi) below).
Our exposition of this Theorem is for the case of
$K = {\CC} \ $ or $\ \RR \ $ \mbox{(e.~g. items ii) and v) ),} but there is
a straightforward algebraic generalization for an arbitrary $K$.
\smallskip

In the Sard-type Theorem below ${\cal G}\ , \ S \ , \ U_{\cal G}$ and bundle
$T_{\cal G} := G|_{\cal G}^{\perp}$ are as constructed above, i.~e. satisfy
the assumptions of Theorem~\ref{Continuation}. Also functions
$L_j \ , \ 1\leq j\leq k \ $, on $U_{\cal G}$ are as constructed above,
i.~e. vanish on $\cal G$ and satisfy (\ref{0}) with $V' = {\cal G}\ $.
\bigskip

\begin{theorem}\label{sard}
({\bf A Sard-type Theorem on singular varieties})
\medskip

For a generic linear
combination $L=\sum _{1\leq j\leq k} c_jL_j$ with  coefficients
$c=(c_1,\dots,c_k)\in {K}^k$
the following properties hold:

i) intersection
$\{L=0\} \cap \Reg_*(\overline S)$ is not empty, dense in $S_{-1}:=\{L=0\}\cap S$
and is smooth  of dimension $\dim(S)-1 \ $;

ii) for any compact (in Euclidean topology on ${K}^n$) set
$C\subset (\overline {S}\cap U_{\cal G})$
and all points $a\in \{L=0\} \cap \Reg_*(\overline S) \cap C$ the norms of
$d(L|_S)(a)=dL(a)|_{T_a(S)}$ are separated from $0$  by
a positive constant (depending on $C$);

iii) the boundary $(\overline{S_{-1}} \setminus S_{-1})\cap U_{\cal G}$
of set $S_{-1}$ in $U_{\cal G}$ coincides with ${\cal G} \ $;

iv) $\Reg(S_{-1})\supset  (S_{-1} \cap \Reg(S))$ and  $S_{-1}$ is G-regular in $U_{\cal G} \ $;

v) for every sequence of points in
$S_{-1}$ converging to a point $a\in \cal G$
such that their tangent
spaces to $S_{-1}$ converge to a subspace $Q$ in the respective
Grassmanian inclusions $T_a({K}^n)\supset Q\supset G_a^{\perp}$ are
valid and therefore also
$$\overline{T(S_{-1})^{\perp}}|_{\cal G}\subset {G}|_{\cal G} \ ;$$

vi) replacing $S_{-1}$ by an irreducible component ${\hat S}_{-1}$ of $S_{-1}$
whose boundary contains $\cal G$ the properties iii)-v) remain valid.
\end{theorem}

\begin{remark}
For the sake of clarity we include though do not make use of the following:

$\bullet$ Of course in ii) of the Lemma above we may equivalently replace "the
norms of $d(L|_S)(a)=dL(a)|_{T_a(S)}$ are separated from $0$" by "the
angles between gradient $\grad L(a)$ of $L$ at $a$ and tangent
spaces $T_a(S)$ to $S$ at $a$ are separated from  $\pi /2 \ $".

$\bullet$ Due to $S$ being irreducible and
$\{ L=0 \} \cap S \neq S$ it follows that
irreducible components of $S_{-1}$ are equidimensional.
\end{remark}

{\bf Deduction of Proposition~\ref{continuation} from Theorem~\ref{sard}}.
We construct sets   ${\hat S}_{-i}:={\hat S}_{-1}({\hat S}_{-i+1})\ ,\ 1\le i\le e:=\dim
(S) -n+k\ $,
consecutively applying $e$ times Theorem~\ref{sard}.
Then due to iii) of  Theorem~\ref{sard}
\begin{eqnarray}\label{1.3}
(\overline{{\hat S}_{-e}}\setminus {\hat S}_{-e})\cap U_{\cal G} \ = \ {\cal G}
\end{eqnarray}
\noindent and,  moreover,
\begin{eqnarray}\label{8}
\overline{T({\hat S}_{-e})^{\perp}}|_{\cal G} = {G}|_{\cal G} \end{eqnarray}
\noindent since  the Gauss map of
${\hat S}_{-e}$ extends (uniquely) as a continuous map to all of $\cal G$
(due to v) of Theorem~\ref{sard}).
Indeed, for every sequence of points from
${\hat S}_{-e}$
converging to a point
$a\in \cal G$
such that their tangent spaces to
${\hat S}_{-e}$
converge to a subspace $Q$ (in the respective Grassmanian), inclusions
$T_a({K}^n)\supset Q\supset G_a^{\perp}$ hold, but  $\dim (Q)=\dim
(G_a^{\perp})=n-k$, and hence  $Q=G_a^{\perp}$.
Therefore due to (\ref{1.3})
${\hat S}_{-e}$
can be enlarged to an irreducible, G-regular and
open in $\overline{{\hat S}_{-e}}$ subset
${\cal G}^+:= {\hat S}_{-e}\cup {\cal G}$
of dimension $n-k$  satisfying (\ref{8}), as required in
Proposition~\ref{continuation}. \bull
\medskip

{\bf Proof of Theorem~\ref{sard}}: {\bf Property vi) follows from
iii)-v)} is straightforward using that $S_{-1}$ is open in its
closure (see Remark~\ref{open}).
\medskip

{\bf We prove iii)} for an arbitrary choice of  $c\in {K}^k$.
Inequalities
$\dim(\overline{(S_{-1})_a})\ge \dim(S)-1\ge n-k>m=\dim({\cal G})$, where
$(S_{-1})_a$ denotes the germ at $a\in \cal G$ of $S_{-1}$ as an analytic set. Using
a similar notation $({\cal G})_a$ for $\cal G$ it follows that
$({\cal G})_a \subset \overline{((\overline{S}\cap \{L=0\}) \setminus {\cal G})_a}$.
On the other hand, $((\overline{S}\cap \{L=0\}) \setminus {\cal G})_a\ =\
((\overline{S} \setminus {\cal G}) \cap \{L=0\})_a\ =\ (S_{-1})_a\ $, since
$(S)_a=(\overline{S}\setminus {\cal G})_a$ due to (\ref{1.2}).
Thus ${\cal G} \subset (\overline{S_{-1}}\cap U_{\cal G})$ and
since also $S\cap {\cal G}=\emptyset \ $,
it follows that $(\overline{S_{-1}} \setminus S_{-1}) \supset \cal G \ $.
Using (\ref{1.2}) it follows  that
${\cal G} \ =\ (\overline{S}\setminus S)\cap U_{\cal G}\ \supset \ (\overline{S_{-1}}
\setminus S_{-1})\cap U_{\cal G}\ \supset \ {\cal G}\ $,
as required in iii).
\medskip

{\bf Properties i) and ii) of Theorem~\ref{sard} imply
both iv) and v).}
Inclusion $\Reg(S_{-1})\supset  \{L=0\} \cap \Reg(S) = (S_{-1} \cap \Reg(S)$ is
a straightforward consequence of i) and ii).
The remainder is a consequence of the following property: if the
limits of two sequences of subspaces of ${K}^n$  exist, then the limit of the respective
intersections of these subspaces also exists and coincides with the
intersection of the limits of the sequences, provided that the
angles between the respective subspaces in the sequences are
separated from $0$ by a positive constant.
\medskip

Thus it remains to prove i) and ii).
\medskip

{\bf Proof of i)}. We have constructed an open in ${K}^n$ set $U_{\cal G}$
and a G-regular irreducible dense subset $S\subset W \cap U_{\cal G}$
of a Lagrangian component of $\{{\cal G}_t\}_{r\leq t\leq k}$ whose
boundary $\overline{S} \setminus S = \cal G$ in $U_{\cal G}$ (see Remark~\ref{open}).
We may assume w.l.o.g. that
$$d(S):=\dim_{K}(\Span\{L_j|_S\}_{1\leq j\leq k})\ge 2 \ ,$$
\noindent where $\Span$ denotes the $K$-linear hull of a family of functions.
Indeed, since $\dim(S)>n-k$ (Corollary~\ref{dense}) it follows
that $d(S)>0 \ $. It remains to exclude the case of $d(S)=1$.
In the latter case we may assume w.l.o.g. that
$\dim(\Span\{L_j|_S\}_{2\le j\le k}) \ge 1$ and then change $L_1$ by adding
to it an appropriate generic element of the square of the ideal $I_{\cal G}$
of all polynomials vanishing on $\cal G$. This would not change the value
of $dL_1$ at the points of $\cal G \ $, but on the other hand $d(S)$ for the new
choice of $L_1$ will increase due to dimension of
$\ I_{\cal G}^2/ I_S\ $ as a vector space over $K$ being infinite, as
required.

We start with an {\it embedded desingularization}
$\sigma:{\cal N}\rightarrow U_{\cal G}$
 of $\overline {S} \cap U_{\cal G} \subset U_{\cal G}$ by means of successive
blowings up along smooth admissible centers \cite{Hironaka},
\cite{BMinvent}, \cite{BMmmj} with 'declared exceptional' hypersurfaces
$H_j,\ 1\leq j \leq k \ $, which we may so declare since the latter are smooth
and they form normal crossings in $U_{\cal G} \ $.
In particular, the following properties hold:

\vspace{2mm}
0. $\sigma: {\cal N} \setminus \sigma^{-1}(\Sing_*(\overline{S})) \to U_{\cal G}  \setminus \Sing_*(\overline{S})$
is an isomorphism;

1. the (so-called) {\it strict transform}
$N:=\overline{\sigma^{-1}((\overline{S}\cap U_{\cal G})
\setminus \sigma (Sing(\sigma)))}$ of $\overline{S}\cap U_{\cal G}$ is smooth;

2. $\Sing_*(\overline{S})=\sigma(\Sing(\sigma))$ and
$\Sing(\sigma)=\sigma^{-1}(\sigma(Sing(\sigma)))=\cup_{i\geq
1}H_{i+k} \ $, where each $H_{i+k}$ is a smooth (so-called) {\it
exceptional} hypersurface and in addition each $H_{i+k}$ is the strict transform
of the set of the critical points of the successive $i$-th
intermediate blowing up;

3. each $H_i\cap N,\ i\geq 1 \ $, is smooth and $\dim(H_i\cap N)=\dim(N)-1$ for $i\geq k+1 \ $;

4. the family $\{H_i\}_{i\geq 0}\ $, where we denote $H_0:=N\ $, forms a normal
crossings in $\cal N\ $.



\vspace{2mm}

For any hypersurface $\{f=0\}\subset U_{\cal G}$ one considers the
{\it strict transform} of $\{f=0\}$
$$\Lambda_{(f)}=\overline{{\sigma}^{-1}(\{f=0\}) \setminus \Sing(\sigma)} \subset {\cal N}$$
\noindent  under  map
$\sigma$.
\begin{remark}\label{strict}
Due to property 2. above the local equation of  $\Lambda_{(f)}$ can be constructed
by factoring out from $f\circ \sigma$ the maximal monomial in exceptional
hypersurfaces. In particular, assume that $f$ depends on parameter $c\in K^k$
and map ${\tilde\sigma} := \sigma \times id \colon {\cal N} \times K^k \rightarrow U_{\cal G} \times K^k \ $.
With $f|_c$ being the evaluation
of $f$ at $c$, hypersurfaces $\Lambda_{(f|_c)}\subset \cal N$ and
$\Lambda_{(f)}\subset {\cal N}\times K^k$ being the
strict transforms under maps $\sigma$ and $\tilde{\sigma}$ respectively, it
follows that if for
a particular value of $c$ hypersurface $\Lambda_{(f)}|_c:=\Lambda_{(f)}\cap
({\cal N}\times \{c\})\subset {\cal N}$ is smooth then
\begin{eqnarray}\label{commute}
\Lambda_{(f|_c)}=\Lambda_{(f)}|_c\ ,
\end{eqnarray}
\noindent where ${\cal N}\times \{c\}$ is identified with $\cal N\ $.
Of course for a sufficiently generic value of $c\in {K}^k$ equality
(\ref{commute}) holds in any case.
\end{remark}

To simplify notation we let $\Lambda_j:=\Lambda_{(L_j)}\subset {\cal
N}\ ,\ 1\le j\le k\ $, and $\Lambda:=\Lambda_{(L)}\subset {\cal
N}\times K^k\ $ (all these hypersurfaces being the strict transforms
under maps $\sigma$ and $\tilde{\sigma}$ respectively).
Hypersurfaces $\Lambda_j \ , \ 1\leq j\leq k \ $, are smooth and
together with $\Sing(\sigma)$ form normal crossings in ${\cal N}$
due to the choice of admissible centers of blowings up (see e.~g.
\cite{BMinvent} or \cite{BMmmj}). In addition, for each $j\ ,\ 1\leq
j\leq k\ ,$ the difference between the divisors of $L_j\circ \sigma$
and $\Lambda_j$ is the exceptional divisor $E_j$ supported on
$Sing(\sigma) = \cup _{i \ge k+1} H_i \subset {\cal N} \ $ (each
divisor being of the form $E_j=\sum _i n_{j,i} [H_i]$ and all
integers $n_{j,i} \geq 0$).

We now, starting with $\cal N \ $, will apply 'combinatorial` blowings up, i.~e.
with centers of all successive blowings up being the intersections of some of
the accumulated exceptional hypersurfaces (possibly including some among
$\Lambda_j \ , \ 1\leq j\leq k \ $).
By means of such blowings up we achieve that the pull back
 of ideal $\cal I$  generated by
$L_j \ , \ 1\leq j\leq k \ $, is principal and, moreover, is
locally generated at any point $a$ by one of the $L_j\circ \sigma \ $,
$\ 1\leq j\leq k$ \cite{BMinvent}. (For
such $j=j(a)$ it follows that $a\not\in \Lambda_j \ $.)
Note that the 'combinatorial part of desingularization' preserves properties 0.-4.
(listed above) of
embedded desingularization of $\overline {S} \cap U_{\cal G} \subset U_{\cal G} \ $.

It follows that $\Lambda$ \noindent is nonsingular. Indeed, for any
point $(x,c)\in \Lambda$ there exists $j \ , \ 1\leq j\leq k \ $,
for which ideal ${\cal I}=(L_j\circ \sigma)$ in a neighbourhood of
point $x\in \cal N \ $. As a consequence, the partial derivative
with respect to $c_j$ of function
$$\lambda:={\sum_{1\leq i\leq k} c_i(L_i\circ \sigma) \over L_j\circ \sigma}$$
\noindent at $(x,c)$ equals 1 and $\{\lambda = 0 \} = \Lambda \ $.

The standard version of Sard's Theorem implies that for a choice of
an appropriate generic $c=(c_1,\dots,c_k)$ the fiber $\Lambda _c$ of
the restriction to $\Lambda$ of the natural  projection  $p\colon
\Lambda \to {K}^k$ is nonsingular in $\sigma^{-1}(U_{\cal G}) \ $.
Note that Sard's Theorem applies because if $x\in {\cal N} \setminus
\Sing(\sigma)$ and $c\neq 0$ then a straightforward calculation
(making use of the linear independence of differentials $dL_j \ , \
\ 1\leq j\leq k \ $, in $U_{\cal G} \ $) shows that the rank of the
Jacobian matrix of the projection $p$ at $(x,c)\in \Lambda$ equals
$k \ $.

To complete the proof of i) we apply Sard's Theorem to the
restriction of $p$ to $\Lambda \cap (N\times {K}^k) \ $.
Note that
$\Lambda \cap (N\times {K}^k) = \{(x,c) \in N\times {K}^k : \lambda (x,c) = 0 \}$
in local coordinates on $N\times {K}^k$ chosen as above and is nonsingular (since
the partial derivative of $\lambda$ with respect to $c_j$ at $(x,c)$ equals $1$ ).
Due to our choice above
$$d(N):=\dim_{K} (\Span (\{L_j\circ \sigma|_{N}\}_{1\leq j\leq k}))=d(S)\ge 2 \ .$$
\noindent Pick $\ L_{j_1}|_S \ , \ \ L_{j_2}|_S \ , \ \ 1\le j_1<j_2\le k \ $,
being linearly independent over $K$. It follows that there is a point
$x\in N\setminus Sing(\sigma)$ and $c_{j_1} \ , \ c_{j_2}\in K$ such that
$$c_{j_1}L_{j_1}(\sigma(x))+c_{j_2}L_{j_2}(\sigma(x))=0 \ ,\quad
c_{j_1}(dL_{j_1})(\sigma(x))+c_{j_2}(dL_{j_2})(\sigma(x))\neq 0$$
\noindent holds. Such point $\ x\in N\setminus Sing(\sigma)$ exists
since otherwise it follows that for all $x\in N\setminus Sing(\sigma)$
$$(L_{j_2} (dL_{j_1})-L_{j_1}(dL_{j_2}))(\sigma(x))=0 \ ,$$
\noindent which would imply a linear dependence of
$L_{j_1}|_S \ , \ \ L_{j_2}|_S$ contrary to their choice.
Set $c_j=0$ for all $j\neq j_1 \ , \ j_2$. Then again  by means of
a straightforward calculation the rank of the Jacobian at $(x,c)$ of
projection $p\colon \Lambda \cap (N\times {K}^k) \to {K}^k$
equals $k$ and therefore Sard's Theorem implies that
$\Lambda _c \cap N$  is nonsingular for appropriate generic $c \ $, where
$N$ is identified with $N\times \{ c \}$. Since $\sigma$ is an isomorphism off $\Sing_*
(\overline{S})$ (which is the property 0. of $\sigma$) it follows that if
$\{L=0\} \cap \Reg_*(\overline S) \neq \emptyset$
then it is a smooth hypersurface of $\Reg_*(\overline S)$ of dimension
$\dim(S)-1 \ $. To complete the proof of i)
it suffices to show that
$\Lambda_c \cap N \not\subset \Sing(\sigma)=\cup_{i\ge 1}H_{i+k}$
and that, moreover,
$\Lambda_c \cap N \setminus \Sing(\sigma)$ is dense in $\Lambda_c \cap N \ $.

Both properties follow by specifying an appropriate generic choice of $c$ further, e.~g.
a choice of $c$ such that $\Lambda_c$ intersects
transversally every $H_J\times \{c\}$ would do, where $H_J=\cap_{j\in J} H_j$ for any
acceptable index set $J \subset \{ i \ge 0 \}\ $. We achieve the latter by
once again applying Sard's Theorem to the restriction of projection
$p$ to $\Lambda \cap (H_J\times {K}^k)$. Of course, for $J$ such that
$p(\Lambda \cap (H_J\times {K}^k))$
is not dense in ${K}^k$ a generic choice of $c\in {K}^k$ implies that
$\Lambda_c \cap (H_J\times {K}^k) =\emptyset \ $, which suffices, and otherwise Sard's
Theorem applies and implies for an appropriate generic choice of $c$ the desired transversality,
which completes the proof of i).

\vspace{3mm}
{\bf Proof of ii)}.  We summarize consequences of application of Sard's Theorem
in the following
\begin{remark}\label{crossing}
For a choice of an appropriate generic $c\in {K}^k$ it follows that the family
$\{ H_i \}_{i \ge 0}$ with $\Lambda_c$ form a normal crossings in ${\cal N}:={\cal N}\times \{c\} \ $.
\end{remark}

For a point $a\in {K}^n$ denote ${\cal L}_a:=\Span(\{\grad L_j(a)
\}_{1\leq j\leq k})\subset {K}^n$. Then ${\cal L}_a+T_a(S)={K}^n$
for all $a\in S$ near any point $b\in \cal G \ $. (Indeed, recall
that $G_b={\cal L}_b^* :=\Span(\{dL_j\}_{1\leq j\leq k}) \ $, due to
(\ref{0}), implying that $k=\dim({\cal L}_b^*) =\dim({\cal L}_a^*)$,
and that $T\supset G_b^{\perp}$ if the limit $T=\lim_{a\rightarrow
b} T_a(S)$ exists, using for the latter inclusion that $S$ is a
dense subset of a Lagrangian component of $\{{\cal G}_t\}_{r\le t\le k}\
$, see Remark~\ref{open}.) Hence $\dim({\cal L}_a \cap
T_a(S))=k+\dim(S)-n$.

There is a natural isomorphism of
$$\Omega _a:={\cal L}_a^*/({\cal L}_a^* \cap T_a(S)^{\perp}) \subset T_a(S)^*$$
\noindent with ${\cal L}_a \cap T_a(S)$ via realization of the functionals on $T_a(S)$ by means of
a scalar product on ${K}^n$. In particular, $\dim (\Omega_a)=k+\dim(S)-n \ $,
$\dim({\cal L}_a^* \cap T_a(S)^{\perp})=n-\dim(S) \ $ and both dimensions do not depend on $a$.

We introduce on ${\cal L}_a^*$ a metric equivalent
to the standard one (over any compact subset of the points $a\in {K}^n$
with $\dim({\cal L}_a^*)=k \ $) by declaring $dL_1,\dots,dL_k$ to be an orthonormal
basis in ${\cal L}_a^* \ $.

For any point ${\tilde b}\in \Lambda _c \cap \Sing(\sigma) \subset
\cal N$ and points ${\tilde a}\in {\cal N}\setminus \Sing(\sigma)$
nearby $\tilde b$ we introduce a metric in $T_{\tilde a}({\cal
N})^*$ as follows. In a neighbourhood of $\tilde b$ the smooth
variety ${\cal N}$ admits a coordinate chart $\cal C$ with the
origin at $\tilde b$ and every exceptional hypersurface $H$
intersecting $\cal C$ by a coordinate hyperplane $\{x_H=0\}$ of
$\cal C$, unless the intersection is empty (one may use here a traditional
complex analytic coordinate chart, or alternatively the notion of an affine
'etale' coordinate chart as in \cite{BMinvent}, \cite{BM99}). In a neighbourhood of
$\tilde b$ the local ideal ${\cal I}_{\tilde b}$ is generated by
a single $L_j \circ \sigma$ for a suitable $j \ $ (as was achieved
by the desingularization above), $\ 1\leq j\leq k\ $, and the function $h:=\lambda|_c$ has
a non-vanishing differential at $\tilde b \ $, since $\Lambda _c \cap N$  is nonsingular
due to the choice of $c$ as shown in the proof of i). We shrink the
neighbourhood $\cal C$ so that $dh$ does not vanish at all points of
$\cal C$. In addition, due to Remarks~\ref{crossing} and \ref{strict}, we may assume
 that $h$ is one of the non-exceptional
coordinates on $\cal C$. We define an auxiliary norm on
$T_{\tilde a}({\cal N})^*\ $ via imposition of the following:
\begin{eqnarray}\label{6}
\left\{ {dx_H \over x_H}, dx_i \right\}_{H,i}\ \mbox{is an orthonormal basis on}\ T_{\tilde a}({\cal N})^*\ ,
\end{eqnarray}
\noindent  where $\{x_H,x_i\}_{H,i}$ are the
coordinates in $\cal C$ with the former ones corresponding to the exceptional hypersurfaces and
the latter $\{x_i\}_i$ being remaining coordinate functions (including function $h\ $).
A straightforward calculation shows that the Hermitian (Riemannian for $K = \RR$) metrics on
${\cal C}\setminus \Sing({\sigma})$ that we have introduced by means of
(\ref{6}) do not depend on the coordinate choices preserving exceptional
hypersurfaces, i.~e. isomorphic for such choices (we do not make use of this
fact), for the case of Hermitian metrics cf. \cite{GM}.

We now will complete the proof of Theorem~\ref{sard} relying on the following lemma, which is
stated in the notations of the preceding paragraph.

\begin{lemma}\label{bound}
The norm of $d(L|_c \circ \sigma)|_{\tilde a} \in T_{\tilde a}({\cal
N})^*$ equals $|L_j\circ \sigma(\tilde a)| \ $, which also majorates the
norm of the linear map $\sigma_{\tilde a}^*\colon \Omega_a \to
T_{\tilde a}({\cal N})^*$ (up to a constant factor depending only on
a choice of $\cal C$) and where ${\tilde a} \in (\Lambda _c \cap {\cal
C})\setminus \Sing(\sigma)$ with $a=\sigma (\tilde a)$.
\end{lemma}

\begin{remark}\label{factor}
The norm of the map  $\sigma_{\tilde a}^*\colon {\cal L}_a^* \to
T_{\tilde a}({\cal N})^*$ equals  the norm of $\sigma_{\tilde
a}^*\colon \Omega_a \to T_{\tilde a}({\cal N})^*$, because the
latter map is the composite of the former one with the natural map
${\cal L}_a^* \to \Omega_a$; therefore it suffices to majorate only
the norm of the former map by $|L_j(a)|$.
\end{remark}

Lemma~\ref{bound} implies a lower bound
depending only on a choice of $\cal C$ on the norms of $(dL)|_S$ at
the points of $\{L=0\}\cap \Reg_*(\overline{S}) \cap \sigma({\cal C})=\Reg_*(\overline{S})\cap
\sigma (\Lambda _c \cap \cal C)$. Since $\sigma$ is a proper map the
item ii) of Theorem~\ref{sard} follows.
 \bull
\bigskip

{\bf Proof of Lemma~\ref{bound}.} As mentioned above $L_j \circ \sigma$ coincides (up to an invertible
function) with $\prod _{{\tilde b}\in H} x_H ^{n_H}$ in $\cal C$ (w.l.o.g. we may assume that
they coincide). Due to Remark~\ref{strict} and using $h({\tilde a})=0$ it follows that
$$d(L|_c \circ \sigma)|_{\tilde a} =d((L_j\circ {\sigma}) \cdot h)|_{\tilde a}=L_j(a) \cdot dh|_{\tilde a} \ .$$
\noindent Due to the choice of the norms on
$T_{\tilde a}({\cal N})^* \ $ (see (\ref{6})), for ${\tilde a} \in {\cal C}\setminus \Sing(\sigma) \ $,
it follows that the norm of $dh|_{\tilde a}$
equals $1\ $. Thus, the norm of $d(L|_c \circ \sigma)|_{\tilde a}$ is
$|L_j(a)| \ $, as required.

It remains to bound the norm of $\sigma_{\tilde a}^*\colon {\cal L}_a^* \to
T_{\tilde a}({\cal N})^*$ (see Remark~\ref{factor}). We observe that the norms of all
$ d(L_i\circ \sigma)|_{\tilde a} \ , \ 1\leq i\leq k \ $, are majorated by $|L_j(a)|$
(up to a constant factor depending only on a choice of $\cal C$)
because $L_j\circ \sigma$ is a common factor of all $L_i\circ \sigma \ , \ 1\leq i\leq k \ $,
in $\cal C$ and the norm of
$ d(L_j\circ \sigma)|_{\tilde a}$ equals $\sqrt{\sum_{{\tilde b}\in H} n_H^2}|L_j(a)| \ $, see
(\ref{6}).
This implies the required upper bound on
the norm of $\sigma_{\tilde a}^*\colon {\cal L}_a^* \to
T_{\tilde a}({\cal N})^* \ , \ $  since the latter is bounded by the maximum of the norms of
the images of the orthonormal basis $\{ dL_i \}_i$ in ${\cal L}_a^* \ $. \bull

\section{Complexity of functorial TWG-stratifications}\label{section6}

One can construct a chain of bundles of vector spaces
$G^{(0)}\subset G^{(1)} \subset \cdots \subset G^{({\rho})}=G$
applying an algorithm for quantifier elimination \cite{CG} to
proceed from $G^{(p)}$ to $G^{(p+1)}\ , \ 0 \le p < {\rho} \ $. This
yields an upper bound $R^{(O(1))}d^{n^{O({\rho})}}$ on complexity
for construction of $G$, where $deg(F)<d$ and $R$ majorates the
bit-size of the coefficients of components $f_i \ , \, 1\leq i\leq l
\ $, of $F = (f_1, \cdot \cdot \cdot , f_l)$ assuming that the
coefficients are, say, algebraic numbers. Note that ${\rho} \le 2n $
(see \cite{Milman03}). Then one can construct quasistrata ${\cal
G}_k$ within the same complexity bound and, if $G$ is Lagrangian, a
functorial TWG-stratification as well (see
Corollary~\ref{strongly}). Note that in an example from
Subsection~\ref{section7.index} the index of stabilization ${\rho}$
grows linearly with $n \ $.

\medskip

\noindent We mention that a similar double-exponential complexity
bound on stratifications (though without properties of universality
nor functoriality) was obtained in \cite{Gabrielov}, \cite{MR},
\cite{Chistov}. On the other hand, there is an obvious exponential
complexity lower bound.
\medskip

It would be interesting to understand, whether this double-exponential bound is
sharp?
\medskip

\section{Examples}\label{section7}

\subsection{A family of ${\bf F: K^N \to K}$ which admit functorial TWG-stratifications}
\label{section7.1}
First we give an example of a family of  polynomials $f\ $, i.~e. $l = 1$ and
$F=(f):K^N \rightarrow K\ $, that admit functorial TWG-stratifications,
which are de facto (in this example) stratifications. (Also,  $G^{(1)} = G \ $,
i.~e. the index of stabilization ${\rho}(f) = 1 \ $.)

Let
$$f=f_n=\sum _{1\leq i\leq j\leq n} A_{i,j}X_iX_j\in K[\{A_{i,j}\},\{X_i\}].$$
Of course $\Sing(f)=\{X_i=0\}_{1\leq i\leq n} \ $. For the sake of brevity let $B$
denote the bundle $G^{(1)}$ of the construction in section 2 that corresponds to
$F := (f) : K^N \to K \ $, where $N = n + {n+1\choose 2} \ $, and $G := G_F \ $.

Any nonsingular $n\times n$ matrix $C$ over $K$ induces an
isomorphism of $K^N \to K^N\ $, which for brevity we also denote $C\
$, and the latter preserves the rank of quadratic forms. Therefore,
for any particular quadratic form $f^{(0)}=\sum _{1\leq i\leq j\leq
n} a^{(0)}_{i,j} X_iX_j$ of a rank $q$ the dimension of the fiber
$B_{f^{(0)}}$ at a point $a^{(0)}=(\{a^{(0)}_{i,j}\},\, \{0\})\in
\Sing(f)$ coincides with the dimension of the fiber $B_{f^{(0)}_q}$
of the quadratic form $f^{(0)}_q=\sum_{1\leq i\leq q}X_i^2 \ $,
e.~g. due to Corollary~\ref{glaeser}.

We identify the set of all quadratic forms of rank $q$ with a
constructible subset ${\cal B}_{k(q)}= (\{a_{i,j}\}, \, \{0\})$ of
$\Sing(f)\ $. A straightforward calculation shows that $\dim({\cal
B}_{k(q)})=qn-q(q-1)/2 \ $. Once again by means of
Corollary~\ref{glaeser} (and of an appropriate isomorphism \mbox{$C
: K^N \to K^N\ $}) it follows that ${\cal B}_{k(q)}$ is smooth and
that  fibers $G_y$ are of the same dimension $k(q)$ at all the
points $y\in {\cal B}_{k(q)}$. (Since $\ l=1\ $ Thom stratification
of $\Sing(F)$ exists by \cite{Hironaka77} and therefore due to (1')
of Lemma~\ref{equivalence} inequality $k(q) \le \codim{\cal
B}_{k(q)}$ holds.) Below we calculate $k(q) \ $, which would allow
us to conclude by making use of Theorem~\ref{lagrangian} that each
${\cal B}_{k(q)}$ is Lagrangian and therefore that $B = G \ $,
${\cal B}_{k(q)} = {\cal G}_{k(q)}$ and that stratification $\{
{\cal B}_{k(q)} \}_{k(q)} \ $, by rank, is a functorial
TWG-stratification.

Consider curves $K \ni t \mapsto K^N$  with $f_q^{(0)}$ at $t=0$ and
defined for any $x^{(0)}\in K^n$ as follows:
\begin{eqnarray}
X_i=t^3x^{(0)}_i \ , \ 1\leq i\leq q \ ; \quad X_j=t^2x^{(0)}_j \ ,
\ q<j\leq n; \quad A_{ii}=1 \ , \ 1\leq i\leq q \ ; \nonumber
\\ \quad A_{jj}=t \ , \ q<j\leq n \ ; \quad  A_{ij}=0 \ , \ i\neq j \nonumber
\end{eqnarray}
A straightforward calculation of the limit along this curve of the
normalized differential $df/||df||$ shows that $\sum_{1\leq i\leq n}
x^{(0)}_idX_i \in B_{f^{(0)}_q} \ $. Consider similarly limits along
curves with the origin at $f_q^{(0)}$ and defined as follows:
$A_{ii}=1 \ , \ 1\leq i\leq q \ $, for all the other pairs of $i \ ,
\ j \ $ with $1\leq i\leq j\leq n$ we set $A_{ij}=t^2 \ $ and also
$X_i=0 \ , \ 1\leq i\leq q$ and $X_j=tx^{(0)}_j \ , \ q<j\leq n \ $.
A straightforward calculation implies that the `coordinate'
projection of $B_{f^{(0)}_q}$ to the subspace spanned by
$\{dA_{ij}\}_{1\leq i\leq j \leq n}$ contains the image under the
degree 2 Veronese map of a point with coordinates $x^{(0)}=(\{0\} \
,\,\{x^{(0)}_j\}_{q<j\leq n}) \in K^n \ $. It follows that subspace
$B_{f^{(0)}_q}$ of $(K^N)^*$ contains vectors $dX_i \ , \ 1\leq
i\leq n \ $, and $dA_{j,s} \ , \ q< j \leq s\leq n \ $, i.~e. $\
k(q) \ge (n+(n-q)(n-q+1)/2) = \codim{\cal B}_{k(q)} \ $, and
therefore $k(q) = \codim{\cal B}_{k(q)} \ $. The latter implying
that each (de facto smooth) quasistratum ${\cal B}_{k(q)}$ is
Lagrangian, $G = B$ and, due to Theorem~\ref{correspondence} and its
Corollary~\ref{strongly}, partition $\{ {\cal B}_{k(q)} \}_{k(q)} \
$, where  $0\leq q\leq n \ $, is the functorial Thom-Whitney-a
stratification of $\Sing(f) \ $. We summarize in the following
\medskip

\begin{proposition}\label{7}
For
$$f=f_n=\sum _{1\leq i\leq j\leq n} A_{i,j}X_iX_j\in K[\{A_{i,j}\},\{X_i\}]$$
the index of stabilization ${\rho}(f) = 1$ and strata ${\cal
B}_{k(q)}=\{(\{a_{ij}\},\{0\}):rk(f)=q\}\subset \Sing(f)$ form a
functorial Thom-Whitney-a stratification with respect to $f$.
\end{proposition}
\smallskip

\subsection{A family of examples of ${\bf F_n: K^{4n+1} \to K}$ with universal \mbox{TWG-stratifications}
and the index of stabilization ${\bf \rho (F_n) = n}$ }\label{section7.index}

Let $q(x,y,u,v,w) := u \cdot x^2 + 2w \cdot x \cdot y + v \cdot y^2$
and produce recursively the following polynomials:

\noindent $q_1:=q(x_1\ ,\ y_1\ ,\ u_1\ ,\ v_1\ ,\ w)\ $,
$\ q_{k+1} := q(x_{k+1}\ ,\ y_{k+1}\ ,\ u_{k+1}\ ,\ v_{k+1}\ ,\ q_k(\cdot))\ ,\ k\ge 1\ $.
Denote
\smallskip

\noindent $\qquad \qquad \qquad \qquad f(\vec{x}\ ,\ \vec{y}\ ,\ \vec{u}\ ,\ \vec{v}\ ,\ w):=q_n(\vec{x}\ ,\ \vec{y}\ ,\ \vec{u}\ ,\ \vec{v}\ ,\ w) \ $,
\smallskip

\noindent where $\vec{x}=(x_1\ ,\ \dots\ ,\ x_n)$
and similarly for
$\vec{y}\ ,\ \vec{u}\ ,\ \vec{v} \ $, i.~e. $f$ depends on $N = 4n +1$ independent variables,
and let $\ h_k:=u_k\cdot v_k -q^2_{k-1}(\cdot)\ , \ 1 \le k \le n \ $.
Then $f=u_n\cdot x^2_n + 2 q_{n-1}\cdot x_n\cdot y_n+v_n \cdot y_n^2\ $ and $\ \Sing(f)=\{x_n=y_n=0\}\ $.
By making use of Corollary~\ref{glaeser} and example from Subsection~\ref{section7.1} it follows that
for points $a\in \Sing(f)$ with $dq_{n-1}(a) \neq 0$ the fibers of bundle $G^{(1)}$ are
\medskip

1. $G_a^{(1)}=\Span\{dx_n\ ;\ dy_n\}\ $ if $\ h_n(a)\neq 0\ $, i.~e.
${\cal G}_2 = \Sing(f)\setminus \{h_n=0\}$ off $\{ dq_{n-1} = 0 \}\ $;
\smallskip

2. $G_a^{(1)}=\Span\{dx_n\ ;\ dy_n\ ;\ dh_n\}\ $ if $h_n(a)=0,\ dh_n(a)\neq 0\ $, i.~e.
off $\{ dq_{n-1} = 0 \}$ quasistratum ${\cal G}_3 =
\Sing(f)\cap \{h_n=0\}\setminus \{ dh_n\neq 0 \}\ $;
\smallskip

3. $G_a^{(1)}=\Span\{dx_n\ ;\ dy_n\ ;\ du_n\ ;\ dv_n\ ;\ dq_{n-1}\}\ $, if $\ h_n(a)=0\ ,\ dh_n(a)=0\ $,
i.~e. ${\cal G}_5 = \Sing(f)\cap \{h_n=0 \ ,\ dh_n=0\}$ off $\{ dq_{n-1} = 0 \}\ $.
\smallskip

4. In the cases 1. and 2. fibers $G_a^{(1)} = (\overline{G^{(0)}})_a \ $, but in the case 3. fibers
$G_a^{(1)}\neq  (\overline{G^{(0)}})_a = $

\noindent $ \left\{ \omega = U_n du_n + V_n dv_n + Q_{n-1} dq_{n-1} +
X_n dx_n + Y_n dy_n\ :\ U_n\cdot V_n=(Q_{n-1}/2)^2\right\} \ ,$
where $\omega$ denotes a $1$-form at $a \ $.
\smallskip

Denote $D_1 := \Span\{dx_n\ ;\ dy_n\ ;\ du_n\ ;\ dv_n\}\ $. Note that
$$df = x_n^2d^{}u_n+y_n^2dv_n+2x_ny_ndq_{n-1}+2(u_nx_n+q_{n-1}y_n)dx_n+2(q_{n-1}x_n+v_ny_n)dy_n\ .$$

Results above rely on elementary calculations of Subsection~\ref{section7.1} summarized below:
\smallskip

\noindent $h_n = \det \left(\begin{array}{c} \ \ u_n \ \ \ q_{n-1}
\\q_{n-1} \ v_n \end{array} \right)\ $ and for any sequence of
points from $K^N$ converging to a point
\smallskip

\noindent $a \in \Sing(f)$ the following holds
\medskip

i) the size of
$\{ {\partial f \over \partial x_n} \ ; \ {\partial f \over \partial y_n} \}$ dominates
$\{ x_n^2\ ,\ y_n^2\ ,\ 2x_n \cdot y_n \}$ at $a$ if $h_n \not \to 0\ $,
\smallskip

ii) the limits of $df/||df||$ are the $1$-forms \mbox{$\omega = U_n
du_n + V_n dv_n + Q_{n-1} dq_{n-1} + X_n dx_n + Y_n dy_n\ $} with $\
U_n\cdot V_n = Q_{n-1}^2/4\ $, since the coefficients of $df\ $ at
$\ du_n\ ,\ dv_n\ ,\ dq_{n-1}$ satisfy $\ x_n^2 \cdot y_n^2 = (2x_n
\cdot y_n)^2/4\ $.
\smallskip

\noindent When $h_n(a) = 0$ the latter also follows from the orthogonality of $\omega \in G_a^{(1)}\ $ to
$\ T_a(\{ h_n = 0 \})$ (see (1') of Lemma~\ref{equivalence}) and
$\ dh_n = v_n \cdot du_n + u_n \cdot dv_n + 2q_{n-1} \cdot dq_{n-1}\ $, implying that
$\omega$ is proportional to $dh_n\ $, while $u_n \cdot v_n = q_{n-1}^2$ for points in
$\{ h_n = 0 \}$.
\smallskip

We now turn to a simple, but crucial observation that the
coefficients of $df\ $ at $\ du_n\ ,\ dv_n\ ,\ dq_{n-1}$ satisfy
inequality $\sqrt{|x_n|^2+|y_n|^2}\ge (\sqrt{2})^{-1} \cdot
|2x_n\cdot y_n|$ and therefore the limits of $df/||df||$ evaluated
at points from $K^N$ that converge to $\Sing(f)\cap \{dq_{n-1}=0\}$
are the $1$-forms with vanishing coefficients at all differentials
of the independent variables on which $q_{n-1}(\cdot)$ depends. In
particular, combining with the preceding summary of the arguments of
Subsection~\ref{section7.1} properties 1. and 2. follow without
making assumption $dq_{n-1}(a) \neq 0$ and also
\smallskip

5. $G_a^{(1)} = D_1\ $ for $\ a\in Z_{n-1} := \Sing (f) \cap \{h_n = 0,\ dh_n= dq_{n-1}=0\}
\subset \{ q_{n-1} = 0\}$
\smallskip

\noindent holds.

{\bf Summarizing}
${\cal G}_2 = \Sing (f) \setminus \{h_n = 0 \}\ ,\ {\cal G}_3 = \Sing (f) \cap \{h_n = 0\ ,
\ dh_n \neq 0 \}$ and with ${\cal G}'_5 := \Sing (f) \cap \{h_n = 0\ ,\ dh_n = 0\ ,
\ dq_{n-1}(a) \neq 0 \}$ bundle $G^{(1)}|_{{\cal G}_2 \cup {\cal G}_3 \cup {\cal G}'_5} =
G|_{{\cal G}_2 \cup {\cal G}_3 \cup {\cal G}'_5}\ $. Also
${\cal G}'_5 = \{ x_n = y_n = u_n = v_n = q_{n-1} = 0\ ,\ dq_{n-1} \neq 0 \}\ ,\ $ and

\noindent $\ Z_{n-1} = \{ x_n = y_n = u_n = v_n = x_{n-1} = y_{n-1} = 0 \} =
\Sing (f) \setminus \left({\cal G}_2 \cup {\cal G}_3  \cup {\cal G}'_5 \right) \ $.
\bigskip

{\bf Detour}. The two Remarks-Examples below are straightforward consequences of the latter
observation and the preceding it summary of the arguments of Subsection~\ref{section7.1}.

\begin{remark}\label{remark.ex}
With notations $\ G = G_{\tilde f}\ ,\ G^{(p)} = G^{(p)}_{\tilde f}\ $ for
a function
$$ {\tilde f} := u \cdot x^2 + 2w^2 \cdot x \cdot y + v \cdot y^2 $$
depending on $5$ variables the following hold:

\noindent inequality $\dim G^{(1)}_a \le 4$ for all $a \in Sing({\tilde f})\ $;
bundles $G\ $ and $\ G^{(1)}$ coincide; quasistrata
${\cal G}_2 = \{ x=y=0\ ,\ u \cdot v - w^4 \neq 0 \}\ ,\
{\cal G}_3 = \{ x=y=0\ ,\ u \cdot v - w^4 = 0\ , \ (u\ ,\ v) \neq 0 \}\ $ and $\ {\cal G}_4 = \{ 0 \}$
are smooth and form \mbox{Thom-Whitney-a} stratification $\cal S$ of $Sing({\tilde f})\ $,
but quasistratum ${\cal G}_4$ is not Lagrangian ($\dim {\cal G}_4 = 0 < 5 - 4\ !$).
Also, $\overline{ G|_{{\cal G}_2}}\ $ and $\overline{ G|_{{\cal G}_3}}\ $ are $\ 5$-dimensional
irreducible components of $G\ $ and $\ G|_{{\cal G}_4}$ is in the closure of $G|_{{\cal G}_3}\ $.
\end{remark}

\begin{remark}\label{ex.sing}
Let non-zero polynomial $g\in K[z_1,\dots,z_m]\ $ and $\ f_g :=
{\tilde f}(u,v,x,y,g(z))\ $, where ${\tilde f}$ is from the
preceding Remark. Denote $G := G_{f_g}\ ,\ G^{(p)} := G^{(p)}_{f_g}\
$. Then for polynomial $f_g$ depending on $m+4$ variables the
following hold:

\noindent inequality $\dim G^{(1)}_a \le 4$ for all $a \in Sing({f_g})\ $;
bundles $G\ $ and $\ G^{(1)}$ coincide; the quasistrata are
$\ {\cal G}_2 = \{ x=y=0\ ,\ u \cdot v - g(z)^4 \neq 0 \},\
{\cal G}_3 = \{ x=y=0\ ,\ u \cdot v - g(z)^4 = 0\ , \ (u,\ v) \neq 0 \}$ and
$\ {\cal G}_4 = \{x=y=u=v=g(z)= 0 \}$; only quasistratum ${\cal G}_4$ is not Lagrangian; the
irreducible components $\overline{ G|_{{\cal G}_2}}\ $ and $\overline{ G|_{{\cal G}_3}}\ $ of $\ G$
are $\ (m+4)$-dimensional and $G|_{{\cal G}_4}$ is in the closure of $G|_{{\cal G}_3}\ $.
{\bf Curiously, an arbitrarily chosen hypersurface $\{g=0\}$ appears as a quasistratum}.
\end{remark}

We now turn to calculation of fibers of $G^{(2)}\ $ for $f\ $. Note that $dq_{n-1} -
2x_{n-1}y_{n-1}dq_{n-2} =$
$$x_{n-1}^2du_{n-1}+y_{n-1}^2dv_{n-1}+2(u_{n-1}x_{n-1}+q_{n-2}y_{n-1})dx_{n-1}+
2(q_{n-2}x_{n-1}+v_{n-1}y_{n-1})dy_{n-1}$$

\noindent and bundles $G=G^{(2)}=G^{(1)}$ off $Z_{n-1} \subset \{ x_{n-1} = y_{n-1} = 0 \}\ $. It follows by making use of Corollary~\ref{glaeser} and of the calculations
like in the summary of the arguments of Subsection~\ref{section7.1} that
 for points $b$
from ${\cal G}'_5$ converging to a point $a \in Z_{n-1} \subset \{ q_{n-1} = 0\ ,\ dq_{n-1} = 0 \}$
with $dq_{n-2} \neq 0$ the span of the limits of the $1$-forms from the fibers $G_b\ $ of $\ G\ $, which includes the
limits of $dq_{n-1}/||dq_{n-1}||\ $, coincides with the fibers of bundle $G^{(2)}\ $, namely:
\smallskip

1'. $G_a^{(2)}=\Span\{dx_{n-1}\ ;\ dy_{n-1}\} \oplus D_1\ $ if $\ h_{n-1}(a)\neq 0\ $, i.~e.
${\cal G}_6 = Z_{n-1}\setminus \{h_{n-1}=0\}$ off $\{ dq_{n-2} = 0 \}\ $;

2'. $G_a^{(2)}=\Span\{dx_{n-1}\ ;\ dy_{n-1}\ ;\ dh_{n-1}\} \oplus D_1\ $
if $h_{n-1}(a)=0,\ dh_{n-1}(a)\neq 0\ $, i.~e.
off $\{ dq_{n-2} = 0 \}$ quasistratum ${\cal G}_7 =
Z_{n-1}\cap \{h_{n-1}=0\}\setminus \{ dh_{n-1}\neq 0 \}\ $;

3'. $G_a^{(2)}=\Span\{dx_{n-1}\ ;\ dy_{n-1}\ ;\ du_{n-1}\ ;\ dv_{n-1}\ ;\
dq_{n-2}\} \oplus D_1\ $, if \mbox{$\ h_{n-1}(a)=0\ $},
\mbox{$dh_{n-1}(a)=0\ ,\ $}
i.~e. ${\cal G}_9 = Z_{n-1}\cap \{h_{n-1}=0 \ ,\ dh_{n-1}=0\}$ off
$\{ dq_{n-2} = 0 \}\ $.

4'. In the cases 1'. and 2'. fibers $G_a^{(2)} = (\overline{G^{(1)}})_a \ $, but in the
case 3'. fibers $G_a^{(2)}\not \subset  (\overline{G^{(1)}})_a$ and the latter
consists of all $1$-forms $\omega\in  G_a^{(2)}$ with coefficients
$U_{n-1}\ ,\ V_{n-1}\ ,\ Q_{n-2}\ $ at $\ du_{n-1}\ ,\ dv_{n-1}\ ,\ dq_{n-2}$
 that satisfy equation $U_{n-1}\cdot V_{n-1}=(Q_{n-2}/2)^2\ $.
Denote $D_2 := \Span\{dx_{n-1}\ ;\ dy_{n-1}\ ;\ du_{n-1}\ ;\ dv_{n-1}\}
\oplus D_1\ $.

Once again, due to the observation that the coefficient of $dq_{n-1}$ at
$dq_{n-2}$ is dominated by its coefficients  at $du_{n-1}\ ,\ dv_{n-1}\ $,
it follows that for points $b\in \Sing(f)$ converging to a point
$a\in \{dq_{n-2}=0\}$ the  limits of the $1$-forms from fibers
$G^{(1)}_b\ $, which include the limits of $dq_{n-1}/||dq_{n-1}||\ $,
consist only of $1$-forms with vanishing coefficients at all differentials
of the independent variables on which $q_{n-2}$ depends. In particular,
properties 1'. and 2'. follow without making assumption $dq_{n-2}(a)\neq 0$
and the fiber of bundle $G^{(2)}\ $ at $\ a$ is
\smallskip

5'. $G^{(2)}_a = D_2\ $ for $\ a\in Z_{n-2}:=Z_{n-1}\cap \{h_{n-1}=0\ ,
\ dh_{n-1}=dq_{n-2}=0\} \subset \{q_{n-2}=0\}\ $.
\smallskip

{\bf Summarizing} $\ {\cal G}_5={\cal G}'_5\ ,
\ {\cal G}_6 = Z_{n-1}\setminus \{h_{n-1}=0\}\ ,
\ {\cal G}_7 = Z_{n-1}\cap \{h_{n-1}=0\ ,\ dh_{n-1}\neq 0 \}$ and with
${\cal G}'_9 := Z_{n-1}\cap \{h_{n-1}=0 \ ,\ dh_{n-1} = 0\ ,\ dq_{n-2} \neq 0 \}\ $ bundle
$G^{(2)}|_{{\cal G}_6\cup {\cal G}_7\cup {\cal G}'_9} =
G|_{{\cal G}_6\cup {\cal G}_7\cup {\cal G}'_9}\ $. Also
${\cal G}'_9 = Z_{n-1}\cap \{u_{n-1} = v_{n-1} = q_{n-2} = 0\ ,\ dq_{n-2} \neq 0 \}\ $, and

\noindent $Z_{n-2} = Z_{n-1}\cap \{u_{n-1} = v_{n-1} = x_{n-2} = y_{n-2} = 0 \} =
Z_{n-1}\setminus \left({\cal G}_6 \cup {\cal G}_7 \cup {\cal G}'_9 \right) \ $.

Thus $G^{(1)} \neq G^{(2)}\ $ and $\ G=G^{(2)}\ $ off $\ Z_{n-2}\ $.
Calculation of fibers of $G^{(p)}\ ,\ p>2$ for points from $Z_{n-2}$
is similar (recursively on $p$), in particular implying that
${\cal G}_9={\cal G}'_9\ $. Summarizing

\begin{proposition}\label{recursion}
Quasistrata $\{{\cal G}_r\}_r$ for polynomial $f$ (in $4n+1$ independent
variables) are smooth, Lagrangian, form a Thom-Whitney-a stratification and
 hence a universal \mbox{TWG-stratification}. The index of stabilization
$\rho(f)\ $ of $\ f$ equals $n\ $.
\end{proposition}

\subsection{Example of ${\bf F: K^5\to K}$ with no universal TWG-stratification}\label{section7.2}

For ${\tilde f}$ from Remark~\ref{remark.ex} we have shown that
there is a non Lagrangian quasistratum and therefore $\Sing({\tilde
f})$ by our main Theorem~\ref{main} does not admit a universal
TWG-stratification. For polynomial ${\tilde f}$ we will reprove this
claim illustrating the proof of Theorem~\ref{main}. In this example
${\cal G}={\cal G}_4\ $, construction of ${\cal G}^+$  is elementary
and we provide it explicitly (cf. Section~\ref{section4.5}). We
choose ${\cal G}^+$ to be a curve defined parametrically by
$\{x=y=0\ ,\ u=v=t^2\ ,\ w=t\}\ $. Then partition of $\Sing({\tilde
f})$ by sets ${\cal B}_2 := {\cal G}_2\ ,\ {\cal B}_3 := {\cal
G}_3\setminus {\cal G}^+\ ,\ {\cal B}_4 := {\cal G}^+$ forms
\mbox{Thom-Whitney-a} stratification $\tilde{\cal S}\ $ with the
associated bundle $B({\tilde{\cal S}})\neq B({\cal S})\ $.

Finally we show that there does not exist a universal
TWG-stratification with respect
to $\tilde f \ $. Assume the contrary, say ${\cal S}^{(0)}$ is
a universal TWG-stratification. Denote by $B({\cal S}^{(0)})$ its bundle
of vector spaces. Proposition~\ref{almost} and
Proposition~\ref{coarser} imply that
$G \subset B({\cal S}^{(0)})\subset (B({\cal S}) \cap B(\tilde {\cal S})) $.
It follows that $G_a=B({\cal S}^{(0)})_a=B({\cal S})_a$ for any point
$a\in {\cal B}_2\cup {\cal B}_3 \ $, while
$G_{0}=B({\cal S}^{(0)})_{0}=B(\tilde {\cal S})_{0}$ is
4-dimensional and is orthogonal to vector
$ {\partial {} \over \partial w} \in T_0(K^5) \ $.
(On the other hand
$B({\cal S})_{0} = (K^5)^{*} \ $). Therefore, ${\cal S}^{(0)}$ being universal
should coincide with $\cal S \ $, but the origin $0$ is not
a Lagrangian stratum of ${\cal S}^{(0)} \ $. Thus our assumption leads to a
contradiction.  Summarizing, we
obtain the following proposition.

\begin{proposition}
There is no universal
\mbox{TWG-stratification} with respect to the polynomial
${\tilde f}=u \cdot x^2 + 2w^2 \cdot x \cdot y + v \cdot y^2 \ $.
\end{proposition}

\subsection{Multiplicities of roots and another functorial TWG-stratification}\label{section7.4}


Let $$f := f_{q+2} = \sum_{0\leq i\leq q} A_i X^iY^{q-i}\in
K[A_0,\dots,A_q,X,Y] \ ,$$ \noindent where $([A_0:\dots
:A_q],X,Y)\in {\PP}^q(K)\times K^2 \ $. In particular, in this
example for every affine chart $\{A_i\neq 0\} \simeq K^q \times K^2
\ , \ 0\le i\le q\ $ of $\ {\PP}^q(K)\times K^2$ we consider mapping $F := f:
K^n \to K\ $, where $n := q+2\ $. Then
$\ \Sing(F)$ admits Thom
stratification and (ii) of Theorem~\ref{lagrangian} applies provided
that all irreducible components of ${\cal G}_k \ , \ n - \dim
(\Sing(F)) \le k \le n$ are of dimension $n - k \ $, which we show
below.

Similarly to the preceding examples $\Sing(f)=\{X=Y=0\}$. Here, in
the original notations of Section~\ref{section2}, we prove for $G :=
G_{f_n}$ (and $G^{(p)} := G^{(p)}_{f_n}$) that index of
stabilization ${\rho}(f_n) = 2 \ $, i.~e. that $G^{(1)} \neq G^{(2)}
= G \ $, bundle $G = G_{f_n}$ is Lagrangian and that $\{{\cal
G}_{k+2}\}_{0\leq k\leq q/2}$ is a universal (and hence functorial)
\mbox{TWG-stratification} with respect to $f_n \ $.

Let us fix a point $a^{(0)}=([a^{(0)}_0:\dots :a^{(0)}_q],0,0) \in
\Sing(f) \ $, for the time being, then  polynomial
\begin{eqnarray}\label{4}
f^{(0)}=\sum_{0\leq i\leq q} a^{(0)}_i X^iY^{q-i}=\prod _j
(b_jX-c_jY)^{m_j}.
\end{eqnarray}
We first verify that for each factor $b_jX-c_jY$ with the multiplicity $m_j\geq 2$ the
fiber of the closure $({\overline{G^{(0)}}})_{a^{(0)}}$ contains
$$v_j := v([c_j:b_j]) = \sum_{0\leq i\leq q} c_j^ib_j^{q-i}dA_i.$$
Consider a line defined (parametrically) as follows:
$$A_i(t)=a^{(0)}_i \ , \ 0\leq i\leq q \ ;\ X(t)=c_jt \ ,\ Y(t)=b_jt \ .$$
Then $\lim _{t\rightarrow 0} df/||df|| $ along this line equals $\
v_j \ $. Conversely, let $v=\sum_{0\leq i\leq q} h_idA_i +cdX+bdY$
with a non-vanishing $(h_0,\dots,h_q)\neq 0$ being the $\lim
_{t\rightarrow 0} df/||df|| $ along a curve
$$(\{A_i(t)\}_{0\leq i\leq q}, X(t), Y(t))\subset {\PP}^q(K)\times K^2$$
with the origin at $ a^{(0)} \ $. Making a suitable
$K$-linear homogeneous transformation $C$ of the 2-dimensional plane
and applying Corollary~\ref{glaeser} we may assume w.l.o.g. that
$ ord_t (X(t)) > ord_t (Y(t)) $ and it suffices to show that
$X^2|f^{(0)} \ . \ $  Assume otherwise, then
$$\ord_t \left\{ {\partial f^{(0)} \over \partial X}, {\partial f^{(0)} \over \partial Y} \right\}=
(q-1)\ord_t (Y(t))<\ord_t (X^iY^{q-i}), \, 0\leq i\leq q \ ,$$ which
contradicts to $ (h_0 \ , \dots \ , h_q) \neq 0 \ $.

Since  vectors $\{v_j\}_j$ form a van-der-Mond matrix and
therefore are linearly independent, it follows

\begin{lemma}\label{step}
For any point $a^{(0)}\in \Sing(f)$  fiber $(G^{(1)})_{a^{(0)}}$ of
bundle $G^{(1)}$ of vector spaces coincides with the linear hull of
vectors $dX,dY$ and $\{v_j\}_j$ for all $j$ with the multiplicity of the factor $b_jX-c_jY$
in $f^{(0)}$ being $m_j\geq 2 \ $ and, moreover, $\dim ((G^{(1)})_{a^{(0)}})-2$
being the number of such $j \ $.
\end{lemma}

For every $v=v([c:b])$ let ${\cal D}^{(l)} (v)$ denote the linear
hull of
$$\left\{ {\partial ^l v \over \partial c^i \partial b^{l-i}} \right\} _{0\leq i \leq l} \ .$$
Then $\{v\}={\cal D}^{(0)}(v)\subset {\cal D}^{(1)}(v)\subset \cdots
\ $ due to the Euler's formula. W.l.o.g. we may assume that $b=1$
(if $b=0$ we exchange the roles of $b$ and $c$) and then  ${\cal
D}^{(l)} (v)$ is the linear hull of the derivatives $\{{\partial ^i
v \over \partial c^i}\}_{0\leq i\leq l} \ $, implying $\dim ({\cal
D}^{(l)} (v))=l+1 \ , \, 0\leq l\leq q \ $.

Below we calculate the limit $\lim _{t\rightarrow 0}
(G^{(1)})_{a^{(t)}} \ $. To that end we consider a curve
$\{a^{(t)}\}_t \subset \Sing(f)$ with the origin at $a^{(0)} \ $,
and assume w.l.o.g. that $a_q^{(t)}=1$ for all $t \ $. Due to
Lemma~\ref{step} we may assume (also w.l.o.g.) that for any $t\neq
0$ the multiplicity of every factor of polynomial
$f^{(t)}=\sum_{0\leq i\leq q} a^{(t)}_i X^iY^{q-i}$ does not exceed
2 and these multiplicities are independent on $t\neq 0 \ $. We may
factorise
$$f^{(t)}=\prod_j \prod_p (X-(c_j+e_{j,p}(t))Y)^{m_{j,p}} \ ,$$
where $1\leq m_{j,p} \leq 2 \ $ and $\ e_{j,p}(t)$ are the appropriate
algebraic functions of $t$ with  $e_{j,p}(0)=0$ for all $j \ , \ p \ $.
Then $\sum_p  m_{j,p}=m_j$ for each $j$ (see (\ref{4})) and we denote
${\overline {m_j}}=\sum_p [m_{j,p}/2] \ $, where by $[m_{j,p}/2]$ we mean
the integral part of $m_{j,p}/2 \ $. Due to  Lemma~\ref{step} it follows that
$\dim ((G^{(1)})_{a^{(t)}})= \sum_j {\overline {m_j}}+2 \ $ for any $\ t\neq 0 $
and  that collection
\begin{eqnarray}\label{2}
\{v([c_j+e_{j,p}(t):1])\}_{m_{j,p}=2} \cup \{dX,\ dY\}
\end{eqnarray}
is a basis of the fiber $(G^{(1)})_{a^{(t)}} \ $.

We claim that
\begin{eqnarray}\label{3}
\lim _{t\rightarrow 0} (G^{(1)})_{a^{(t)}}=\bigoplus_j {\cal
D}^{({\overline {m_j}}-1)}(v([c_j:1])) \oplus \Span\{dX,\ dY\} \ .
\end{eqnarray}
To that end we observe that the right-hand side of (\ref{3}) is indeed the direct
sum of the vector spaces due to the Hermite's interpolation  (which
interpolates uniquely a polynomial in terms of the values of its several
consecutive derivatives at the given points, cf. Appendix).
Therefore the dimension of the right-hand side equals
$\ \sum_j {\overline {m_j}}+2 \ $ and to complete the proof of (\ref{3}) it suffices
to verify that the left-hand side of (\ref{3}) contains its right-hand side.

To this end fix $j \ $, denote $m := {\overline {m_j}} \ $ and let
$$E^{(i)} := (\{e^i_{j,p}(t)\}_{1\leq p\leq m})^T \in K^m \ , \ i\geq 0 \ ,$$
where all $p$ satisfy $m_{j,p}=2$ (see (\ref{2})). Let $E$ be the $m\times m$ van-der-Mond
matrix with the columns $E^{(i)} \ , \  0\leq i\leq m-1 \ $. Consider an arbitrary
$w=(w_0 \ , \ \dots \ , \ w_{m-1}) \in K^m$ and let
$u := (\{u_p\}_{1\leq p\leq m}):=wE^{-1} \ $.
Since $E^{-1}E^{(i)}(0)=0$ for every $i\geq m$ it follows for
$u^{(i)}(t) := u \cdot E^{(i)}(t) \ $  that  $\ u^{(i)}(0)=0 \ $. Therefore
$$\sum_{1\leq p\leq m} u_pv([c_j+e_{j,p}(t):1])=\sum_{0\leq s\leq m-1} {w_s \over s!}
{d^sv([c_j:1]) \over dc^s} + \sum_{m\leq i\leq q} {u^{(i)} \over
i!}{d^iv([c_j:1]) \over dc^i} \ .$$ Claim (\ref{3}) follows by
letting $t=0$ in the right-hand side of the latter (in view of the
choice of $w$ as an 'arbitrary' in $K^m$).

We now specify the choice of curve $\{a^{(t)}\}_t \ $  so that for every $j$
equality $\ {\overline {m_j}}=[m_j/2] \ $ holds
(see (\ref{4})), in other words  $m_{j,p}=2 \ $ for $\ {\overline {m_j}} \ $
number of $\ p$'s and, moreover, in the case when number $m_j$ is odd that
$m_{j,p_0} = 1$ for a single $p_0 \ $. Then due to (\ref{3}) it follows

\begin{proposition}\label{closed}
For any point $a^{(0)}\in \Sing(f)$ the fiber
$$(G^{(1)})_{a^{(0)}}=\bigoplus_j {\cal D}^{([m_j/2]-1)}(v([c_j:1])) \oplus \Span\{dX,\ dY\}$$
is a vector space of the dimension $\sum_j [m_j/2]+2$ (see
(\ref{4})). In particular, bundle \mbox{$G := G_f = {\overline {G^{(1)}}} \ $.}
\end{proposition}

Finally, we establish that $G$ is Lagrangian. For every $ k \ , \
0\leq k\leq q/2 \ ,$ let
$${\cal G}_{k+2}^{(0)}:=\{a^{(0)}\in \Sing(f): f^{(0)}=\prod_{1\leq j\leq k} (X-c_jY)^2
\cdot \prod_{k<s\leq q-k} (X-c_sY)\} \ ,$$ i.~e. $f^{(0)}$ has $k$
factors of multiplicity $2$ and $q-2k$ factors of multiplicity $1 \
$. Proposition~\ref{closed} implies that ${\cal G}_{k+2}^{(0)}
\subset {\cal G}_{k+2}$ (see Definition~\ref{definition}) and,
moreover, that ${\cal G}_{k+2}^{(0)}$ is dense in ${\cal G}_{k+2}$.
On the other hand, ${\cal G}_{k+2}^{(0)}$ is open and is isomorphic
to the set of all orbits of the group $\Sym(k)\times \Sym(q-2k)$
acting on a set
$${\cal Z}:= K^{q-k}\setminus (\bigcup _{1\leq i<j\leq q-k} \{Z_i=Z_j\}) \ ,$$
where $\Sym(k)$ permutes the first $k$ coordinates $Z_1\ , \ \dots \
,\ Z_k \ $ and $\ \Sym(q-2k)$ permutes the last $q-2k$ coordinates
$Z_{k+1}\ , \ \dots \ ,\ Z_{q-k} \ $. It follows $\dim({\cal
G}_{k+2}^{(0)}) =q-k \ $. Moreover, ${\cal G}_{k+2}^{(0)}=H({\cal
Z}) \ $, where $H$ maps $Z_1\ ,\ \dots \ ,\ Z_k$ to double roots of
$\ f^{(0)} \ $ and $Z_{k+1}\ ,\ \dots \ ,\ Z_{q-k}$ to single roots.
It follows that ${\cal G}_{k+2}^{(0)}$ is irreducible. Finally,
since in this example $\Sing(F)$ admits Thom stratification,
quasistrata ${\cal G}_{k+2}$ are irreducible and of dimension
$n-k-2$ item (ii) of Theorem~\ref{lagrangian} and hence
Corollary~\ref{strongly} apply
and imply the following

\begin{theorem}
Index of stabilization ${\rho}(f_{q+2}) = 2 \ $, bundle $G =
G_{f_{q+2}}$ is Lagrangian and $\{{\cal G}_{k+2}\}_{0\leq k\leq
q/2}$ is a functorial \mbox{TWG-stratification} with respect to
$f_{q+2} \ $.
\end{theorem}

\section{Appendix. Complexity of extension to a Gauss regular subvariety
with a prescribed tangent bundle over singularities}

Here we estimate complexity of an algorithm of extending of a (smooth) singular
locus of an algebraic variety to a Gauss regular subvariety with a prescribed
tangent bundle over the singularities of the variety (see Section~\ref{section4.5}).
We follow the notations of Sections~\ref{section4}, \ref{section4.5}
with an exception that we use $K$ rather than $\CC \ $.
The input for this algorithm is a family of polynomials
$g_p\ ,\ M_{j+m,i+m} \in K_0[X_1,\dots,X_n] \ $
with $\ p\ge 0\ ,\ i\ ,\ j$ for a subfield $K_0\subset K \ $. For the sake of complexity
bounds we assume that elements of $K_0$ can be represented algorithmically, e.~g.
one may use here the field of rational or algebraic numbers in place of $K_0 \ $, cf. \cite{CG}.
We assume the following representation for an algebraic variety
$S=\{g_0 \cdot g_1\neq 0 \ ,\ g_p=0\}_{p\ge 2}$ and its (smooth) singular locus
${\cal G}=\{g_0\neq 0,\ g_p=0\}_{p\ge 1} \ $, which also is its boundary in $\{ g_0 \neq 0 \}$
(see Remark~\ref{open}). The output of the algorithm is a Gauss regular subvariety
${\cal G}^+$ of $\overline {S} \cap \{ g_0 \neq 0 \}$ (see Proposition~\ref{continuation}).

Basically the algorithm consists of 3 subroutines. The first one is choosing a Noether
normalisation $\pi$ for ${\cal G} \ $. The second one is an {\it implicit parametric
interpolation} of polynomials $L_j$ from Section~\ref{section4.5}. (We refer to the
latter as implicit because the interpolation data are given over the subsets of
points from $\cal G$ and thus the data appear implicitly.)
The third subroutine is a construction of ${\cal G}^+$ proper. To this end we may exploit
a choice of algebraically independent coefficients $c_1,\dots,c_k$ at each consecutive
application of Theorem~\ref{sard} and thereafter to construct an irreducible component
containing $\cal G$ of the resulting intersection with
$\overline {S} \cap \{ g_0 \neq 0 \}$ (cf.  vi) of Theorem~\ref{sard} and the deduction of
Proposition~\ref{continuation}).
Complexity bounds for Noether normalisation and for constructing irreducible components
one may find in  \cite{Logar}, and in \cite{CG} respectively. We observe that the third
subroutine depends only on the complexity of finding irreducible components. We
therefore focus on an algorithm for a parametric interpolation. In fact,
we design an algorithm for interpolation over the parameters varying in $K^m \ $, whereas
for the purposes of Section~\ref{section4.5} it suffices to have the parameters varying in
an open subset ${\cal U}'\subset K^m \ $, which would have simplified the algorithm.

To formulate the complexity bounds we assume that $\deg(g_p)<\delta \ ,\ \deg(M_{j+m,i+m})<\Delta \ $
for all $\ p\ ,\ i \ ,\ j$ and the total number of bits in representation of the coefficients
(in  $K_0$)
of polynomials $g_p\ ,\ M_{j+m,i+m}$ does not exceed $R$. Our main result here is the
following

\begin{proposition}\label{interpolation}
One can interpolate polynomials $L_j$ as required in Section~\ref{section4.5} and, moreover,
under assumptions listed in the preceding paragraph $\deg(L_j)<\Delta\delta^{O(n)}$
is a bound on the degrees of the resulting $L_j \ $. Complexity bound for this
interpolation algorithm is $(R\Delta^n\delta^{n^2})^{O(1)}$.
\end{proposition}

Combining with the complexity bounds for the first and the third subroutines
it follows

\begin{corollary}
The complexity of the algorithm constructing ${\cal G}^+$ is
bounded by $$R^{O(1)}(\Delta \delta)^{n^{O(1)}}.$$
\end{corollary}

{\bf Proof of Proposition~\ref{interpolation}.}
We first consider a {\it non-parametrical} interpolation.

\begin{lemma}\label{non-parametric}
Let $v_1,\dots,v_t\in K^{n-m} \ $ and $\ w_q^{(i)}\in K \ ,\ 1\leq q\leq t\ ,\ 0\leq i\leq n-m \ $ .
There exists a polynomial $A\in K[X_{m+1},\dots,X_n] \ $ of $\ \deg(A)<2t(n-m)$ such that
$$A(v_q)=w_q^{(0)}\ ,\ {\partial A\over \partial X_{i+m}}(v_q)=w_q^{(i)}\ ,\ 1\leq q\leq t\ ,\ 1\leq i\leq n-m \ .$$
\end{lemma}

{\bf Proof.} By making an appropriate linear change of the coordinates in $K^m$ we may assume w.l.o.g. that
$\ v_{q_1}^{(i)}\neq v_{q_2}^{(i)} \ , \ 1\leq q_1<q_2\leq t \ ,\ 1\leq i\leq n-m \ $,
where $v_q=(v_q^{(1)},\dots,v_q^{(n-m)}) \ ,\  1\leq q\leq t \ $. Consider a polynomial
$$A_{q_0}=\prod_{q\neq q_0,1\leq i\leq n-m} (X_{i+m}-v_q^{(i)})^2 \cdot \left(\sum_{1\leq i\leq n-m} a_i(X_{i+m}-
v_{q_0}^{(i)})+a_0\right) \ ,\  1\leq q_0 \leq t$$ \noindent with
indeterminate coefficients $a_i \ ,\ 0\leq i\leq n-m \ $. Then
$A_{q_0}(v_q) = {\partial A_{q_0} \over \partial X_{i+m}}(v_q)=0 \
,\ 1\leq i\leq n-m \ $, for every $q\neq q_0 \ $. Equation
$A_{q_0}(v_{q_0}) = w_{q_0}^{(0)}$ uniquely determines $a_0 \ $.
Furthermore equation ${\partial A_{q_0} \over \partial
X_{i+m}}(v_{q_0})=w_{q_0}^{(i)}$ uniquely determines $a_i \ ,\ 1\leq
i\leq n-m \ $. Finally we let $A := \sum_{1\leq q\leq t} A_q$. \bull
\bigskip

Of course one can in the same vain interpolate the higher derivatives as well.

We now consider a {\it parametric} interpolation.
Due to B\'ezout inequality $\deg(\overline{\cal G}) < \delta^n\ $, we introduce a polynomial
$${\cal A}=\sum_{0\leq e_1+\cdots+e_{n-m}\leq 2(n-m)\delta^n} A_EX_{m+1}^{e_1}\cdots X_n^{e_{n-m}}$$
\noindent with indeterminate coefficients $a:=\{A_E\}_E \ ,\ E=(e_1,\dots,e_{n-m})$ and a
quantifier-free formula $\Phi(u,v,a)$ of the theory of algebraically closed fields which says that
$$\mbox{if}\ v\in {\cal G} \ ,\ \pi(v)=u\in K^m\ \mbox{then}\ {\cal A}(v)=0 \ ,
\ {\partial {\cal A} \over \partial X_{i+m}}(v)=M_{j+m,i+m}(v) \ ,\ 1\leq i\leq n-m$$
\noindent for some $j \ , \ 1\leq j\leq k \ $ (we fix $j$ for the time being). Then the formula
$\ \forall u\ \exists a\ \forall v\ \Phi \ $
is true due to Lemma~\ref{non-parametric}.

An algorithm from \cite{GV} yields a representation of
$\pi^{-1}(u)\cap {\cal G}$ commonly refered to as a ``shape lemma''. Applied
to a system $\{g_p=0, g_0\neq 0\}_{p>0}$ the output of this algorithm
is a partition of $K^m=\cup_{\beta}U_{\beta}$ into constructible subsets such that for each
$\beta$ there are a linear combination $\alpha=\sum_{1\leq i\leq n-m}\alpha_{i,\beta}v^{(i)}$
of coordinates $v^{(i)} \ ,\ 1\leq i\leq n-m \ $, with integer coefficients $\alpha_{i,\beta}$
and rational functions $\phi \ ,\ \phi_i\in K_0(X_1,\dots,X_m)[Y] \ ,\ 1\leq i\leq n-m \ $,
for which the following holds:

$\bullet$ for any  $u\in U_{\beta}$ and any
$v = (u,v^{(1)}\ ,\ \dots \ ,\ v^{(n-m)}) \in \pi^{-1}(u)\cap {\cal G} \ $
equalities $v^{(i_0)} = \phi_{i_0}(u \ ,\ \alpha) \ ,\ 1\leq i_0\leq n-m \ $, take place,
i.~e. $\alpha$ is a primitive element of the field $K_0(u\ ,\ v^{(1)}\ ,\ \dots \ ,\ v^{(n-m)})$
over $K_0(u) \ $;

$\bullet$  the roots of a univariate polynomial $\phi(u,Y)$ are exactly the values of $\alpha$
while ranging over points $v\in \pi^{-1}(u)\cap {\cal G} \ $.

Furthermore, in formula $\Phi$ we replace $v^{(i_0)} \ ,\ 1\leq
i_0\leq n-m \ $, by $\phi_{i_0}(u,\alpha)$ and divide the resulting
polynomials ${\cal A}(\alpha)$ and $\left({\partial {\cal A}\over
\partial X_{i+m}}(\alpha)-M_{j+m,i+m}(\alpha) \right)$ by polynomial
$\phi(u,\alpha) \ $ (with the remainders as polynomials in
$\alpha$). Then system $\Phi_1$ obtained by equating to zero all
coefficients of the remainders at the powers of $\alpha$ is
equivalent to formula $\forall v\ \Phi \ $, for any $u\in U_{\beta}
\ $.

One may consider $\Phi_1$ as a linear system with respect to
variables $\ a\ $ and apply to $\Phi_1$ an algorithm of {\it
parametric Gaussian elimination} (see e.~g.  \cite{GV}). It yields a
refinement $K^m=\cup_{\beta'}U_{\beta'}'$ of partition
$\cup_{\beta}U_{\beta}$ into constructible subsets such that for each
$\beta'$ and for every multiindex $E$ there is rational function
$a_E\in K_0(X_1,\dots,X_m)$ such that for any $u\in U_{\beta'}'$ the
array of coefficients $a(u)=\{a_E(u)\}_E$ fulfils $\Phi_1 \ $. For a
choice of the unique $\beta'$ for which $U_{\beta'}'$ is dense in
$K^m $ the rational function
$$L_j=\sum_{0\leq e_1+\cdots+e_{n-m}\leq 2(n-m)\delta^n} a_EX_{m+1}^{e_1}\cdots X_n^{e_{n-m}}$$
\noindent corresponding to this $\beta'$ is as required in Section~\ref{section4.5}.

Finally we address the complexity issue.
In the construction of the ``shape lemma'' above $\deg(\phi)\ ,\ \deg(\phi_i)$ are bounded by
$\delta^{O(n)}$ as well as the degrees of the polynomials representing $\{U_{\beta}\}_{\beta} \ $,
while the number of  $\{U_{\beta}\} \ $, the total sum of sizes of the coefficients of these
polynomials and the complexity of the algorithm do not exceed $R^{O(1)}\delta^{O(n^2)} \ $ \cite{GV}.
Therefore the degrees of the polynomials occuring in $\Phi_1$ are bounded by $\Delta\delta^{O(n)} \ $,
while the number of the polynomials, the total sum of sizes of their coefficients and the
complexity of constructing $\Phi_1$ do not exceed $(R\Delta^n\delta^{n^2})^{O(1)} \ $.
At the stage of
applying the parametric Gaussian elimination to $\Phi_1$ the bounds are similar.
Proposition is proved. \bull

{\bf Acknowledgements.} The authors are thankful to the Max-Planck
Institut f\"ur Mathematik, Bonn for its hospitality during writing
the paper.


\begin{thebibliography}{99}
\bibitem{BMinvent}
E.~Bierstone, P.~Milman, {\em Canonical desingularization in characteristic zero by blowing up
the maximum strata of a local invariant}, Invent. Math., {\bf 128} (1997), 207--302.
\bibitem{BM99}
E.~Bierstone, P.~Milman, {\em Standard Basis along a Samuel stratum, and implicit differentiation},
Fields Inst. Commun., {\bf 24}, Amer. Math. Soc., Providence, RI, (1999), 81--113.
\bibitem{BMmmj}
E.~Bierstone, P.~Milman, {\em Desingularization algorithms. I. Role of exceptional divisors},
Mosc. Math. J., {\bf 3} (2003), 751--805.
\bibitem{Milman03}
E.~Bierstone, P.~Milman, W.~Paw\l ucki, {\em Differential functions defined in closed
sets. A problem of Whitney}, Invent. Math., {\bf 151} (2003), 329--352.
\bibitem{Chistov}
A.~Chistov, {\em Efficient smooth stratification of an algebraic variety of characteristic
zero and its applications}, J. Math. Sci., {\bf 113} (2003), 689--717.
\bibitem{Gabrielov}
A.~Gabrielov, N.~Vorobjov, {\em Complexity of stratifications of semi-Pfaffian sets},
Discrete Comput. Geom., {\bf 14} (1995), 71--91.
\bibitem{Loo}
C.~G.~Gibson, K.~Wirthm\"uller, A.~A.~du Plessis, E.~J.~N.~Looijenga, {\em Topological
stability of smooth mappings}, Lect. Notes Math., {\bf 552}, 1976.
\bibitem{Glaeser}
G.~Glaeser, {\em Etude de quelques alg\`ebres tayloriennes}, J. Analyse Math.,
{\bf 6} (1958), 1--124.
\bibitem{Mac}
M.~Goresky, R.~MacPherson, {\em Stratified Morse theory}, Springer,
1988.
\bibitem{GM}
C.~Grant, P.~Milman, {\em Metrics for singular analytic spaces},
Pacific J. Math. {\bf 168}, (1995), 61--156.
\bibitem{CG}
D.~Grigoriev, {\em Computational complexity in polynomial algebra}, Proc. of
International Congress of Mathematicians,  Berkeley, vol.~2, (1986), 1452--1460.
\bibitem{GV}
D.~Grigoriev, N.~Vorobjov, {\em Bounds on numbers of vectors of
multiplicities for polynomials which are easy to compute},  Proc.
ACM Intern. Conf. Symbolic and Algebraic Computations, Scotland,
(2000), 137--145.
\bibitem{Merle}
J.~P.~Henry, M.~Merle, C.~Sabbah, {\em Sur la condition de Thom stricte pour
un morphisme analytique complexe}, Ann. sci. de l'Ecole Normale Sup\'er.
S\'er. 4, {\bf 17} (1984), 227--268.
\bibitem{Hironaka}
H.~Hironaka, {\em Resolution of singularities of an algebraic variety over a field of
characteristic zero I, II}, Ann. Math., {\bf 79} (1964), 109--326.
\bibitem{Hironaka73}
H.~Hironaka, {\em Number Theory, Algebraic Geometry and Commutative Algebra}, Volume in honour
of Y.~Akizuki, publ. Kinokuniya, Tokyo, 1973.
\bibitem{Hironaka77}
H.~Hironaka, {\em Stratification and flatness. Real and complex singularities}, Proc.
Nordic Summer School, Oslo, publ. Sijhoff \& Noordhoff, (1977), 199--265.
\bibitem{Kaloshin}
V.~Kaloshin, {\em Around the Hilbert-Arnold problem}, CRM Monogr. Ser., {\bf 24} (2005),
111--162.
\bibitem{Kuo}
T.~C.~Kuo, {\em The ratio test for analytic Whitney stratifications},
Lecture Notes in Mathematics, {\bf 192}, (1969/70), 141--149.
\bibitem{Logar}
A.~Logar, {\em A computational proof of the Noether normalization lemma}, Lect. Notes
Comput. Sci., {\bf 357}, (1988), 259--273.
\bibitem{MR}
T.~Mostowski, E.~Rannou, {\em Complexity of the computation of the canonical
Whitney stratification of an algebraic set in $\CC^n$}, Lect. Notes Comput. Sci.,
{\bf 539} (1991), 281--291.
\bibitem{Adam}
A.~Parusinski, {\em On the bifurcation set of complex polynomial with isolated singularities
at infinity}, Compositio Mathematica, {\bf 97}, (1995), 369-384.
\bibitem{Thom}
R.~Thom, {\em Propri\'et\'es Diff\'erentielles Locales des Ensembles Analytiques},
S\'emin. Bourbaki, {\bf 281} (1964/5).
\bibitem{Wall}
C.~T.~C.~Wall, {\em Regular Stratifications}, Lect. Notes Math., {\bf 468} (1974), 332--344.
\bibitem{Whitney}
H.~Whitney, {\em Tangents to an analytic variety}, Ann. Math., {\bf 81} (1965),
496--549.
\end{thebibliography}
\end{document}